\documentclass{article}

\usepackage{amssymb,amsmath,amsthm}
\usepackage[all]{xy} 
\usepackage{xspace}
\usepackage{pxfonts} 
\newcommand{\boxEff}{\makebox[170pt]{Cartesian effect categories}}
\newcommand{\boxArr}{\makebox[230pt]{Arrows}}

\newcommand{\bC}{\mathbf{C}}
\newcommand{\bV}{\mathbf{V}}
\newcommand{\Set}{\mathbf{Set}}
\newcommand{\Part}{\mathbf{Part}}

\newcommand{\cD}{\mathcal{D}}

\newcommand{\eqst}{\equiv}
\newcommand{\eqwe}{\lesssim}
\newcommand{\weeq}{\gtrsim}

\newcommand{\id}{\mathrm{id}}
\newcommand{\app}{\mathrm{app}}
	
\newcommand{\parto}{\rightharpoonup}

\newcommand{\rupto}[1]{\stackrel{#1}{\longrightarrow}}
\newcommand{\lupto}[1]{\stackrel{#1}{\longleftarrow}}
\newcommand{\cone}[5]{#1\lupto{#2}#3\rupto{#4}#5}

\newcommand{\srupto}[1]{\stackrel{#1}{\rightsquigarrow}}
\newcommand{\slupto}[1]{\stackrel{#1}{\leftsquigarrow}} 
\newcommand{\scones}[5]{#1\slupto{#2}#3\srupto{#4}#5}
\newcommand{\scone}[5]{#1\slupto{#2}#3\rupto{#4}#5}
\newcommand{\cones}[5]{#1\lupto{#2}#3\srupto{#4}#5}

\newcommand{\sto}{\rightsquigarrow}

\newcommand{\tuple}[1]{\langle #1 \rangle}
\newcommand{\np}{\tuple{\,}}

\newcommand{\arr}{\mathtt{arr}}
\newcommand{\first}{\mathtt{first}}
\newcommand{\scond}{\mathtt{second}}
\newcommand{\acomp}{>\!\!>\!\!>}
\newcommand{\altimes}{*\!\!*\!\!*}
\newcommand{\atuple}{\&\!\!\&\!\!\&}
\newcommand{\assoc}{\mathtt{assoc}}
\newcommand{\swap}{\mathtt{swap}}
\newcommand{\fst}{\mathtt{fst}}

\newcommand{\tA}{\mathtt{A}}

\newcommand{\tr}{\mathit{trans}}
\newcommand{\sym}{\mathit{sym}}
\newcommand{\subst}{\mathit{subst}}
\newcommand{\repl}{\mathit{repl}}
\newcommand{\comp}{\mathit{comp}}

\begin{document}

\title{Sequential products in effect categories}
\newcommand{\email}[1]{{\small\texttt #1}}
\author{Jean-Guillaume Dumas \\ LJK, University of Grenoble, France.
\email{Jean-Guillaume.Dumas@imag.fr}
\\[2mm] Dominique Duval \\ LJK, University of Grenoble, France.
\email{Dominique.Duval@imag.fr}
\\[2mm] Jean-Claude Reynaud \\ Malhiver, 38640 Claix, France.
\email{Jean-Claude.Reynaud@imag.fr} } 

\date{July 4., 2007}


\theoremstyle{plain}

\newtheorem{thm}{Theorem}[section]
\newtheorem{cor}[thm]{Corollary}
\newtheorem{lem}[thm]{Lemma}
\newtheorem{prop}[thm]{Proposition}
\newtheorem{asm}[thm]{Assumption}

\theoremstyle{definition}

\newtheorem{rem}[thm]{Remark}
\newtheorem{rems}[thm]{Remarks}
\newtheorem{exa}[thm]{Example}
\newtheorem{exas}[thm]{Examples}
\newtheorem{defi}[thm]{Definition}
\newtheorem{conv}[thm]{Convention}
\newtheorem{conj}[thm]{Conjecture}
\newtheorem{prob}[thm]{Problem}
\newtheorem{oprob}[thm]{Open Problem}
\newtheorem{algo}[thm]{Algorithm}
\newtheorem{obs}[thm]{Observation}
\newtheorem{qu}[thm]{Question}

\numberwithin{equation}{section}



\maketitle

\begin{abstract}
\noindent A new categorical framework is provided for 
dealing with multiple arguments in a programming language
with effects, for example in a language with imperative features.
Like related frameworks (Monads, Arrows, Freyd categories), 
we distinguish two kinds of functions.
In addition, we also distinguish two kinds of equations. 
Then, we are able to define a kind of product, 
that generalizes the usual categorical product. 
This yields a powerful tool for deriving many results
about languages with effects. 
\end{abstract}

\section{Introduction}


The aim of this paper is to provide a new categorical framework
dealing with multiple arguments in a programming language
with effects, for example in a language with imperative features.
In our \emph{cartesian effect categories},  
as in other related frameworks (Monads, Arrows, Freyd categories), 
two kinds of functions are distinguished. 
The new feature here is that 
two kinds of equations are also distinguished.
Then, we define a kind of product, that is mapped to 
the usual categorical product when the distinctions (between functions
and between equations) are forgotten.
In addition, 
we prove that cartesian effect categories determine Arrows. 

\medskip 

A well-established framework for dealing with computational effects
is the notion of \emph{strong monads}, that is used in Haskell 
\cite{Moggi91,Wadler93}.
Monads have been generalized 
on the categorical side to \emph{Freyd categories} \cite{PowerRobinson97}
and on the functional programming side to \emph{Arrows} \cite{Hughes00}. 
The claims that Arrows generalize Monads and 
that Arrows are Freyd categories are made precise in \cite{HeunenJacobs06}.
In all these frameworks, effect-free functions are distinguished
among all functions, 
generalizing the distinction of \emph{values} among all \emph{computations}
in \cite{Moggi91}.
In this paper, as in \cite{BentonHyland03,HeunenJacobs06}, effect-free functions
are called \emph{pure} functions;
however, the symbols $\bC$ and $\bV$, 
that are used for the category of all functions 
and for the subcategory of pure functions, respectively, 
are reminiscent of Moggi's terminology. 

In all these frameworks, one major issue 
is about the order of evaluation of the arguments of multivariate operations. 
When there is no effect, the order does not matter, 
and the notion of product in a cartesian category provides 
a relevant framework. 
So, the category $\bV$ is cartesian, 
and \emph{products of pure funtions} are defined by 
the usual characteristic property of products. 
But, when effects do occur, the order of evaluation of the arguments
becomes fundamental, 
which cannot be dealt with the categorical product.
So, the category $\bC$ is not cartesian,
and products of functions do not make sense, in general. 
However, some kind of \emph{sequential product of computations}
should make sense, 
in order to evaluate the arguments in a given order.
This is usually defined, by composition, from 
some kinds of products of a computation with an identity.
This is performed 
by the \emph{strength} of the monad \cite{Moggi91}, 
by the \emph{symmetric premonoidal category} of the Freyd category 
\cite{PowerRobinson97},
and by the \emph{first} operator of Arrows \cite{Hughes00}.

\medskip

In this paper, the framework of \emph{cartesian effect categories} 
is introduced. 
We still distinguish two kinds of functions: 
pure functions among arbitrary functions,
that form two categories $\bV$ and $\bC$, 
with $\bV$ a subcategory of $\bC$, and $\bV$ cartesian.
Let us say that the functions are \emph{decorated},
either as pure or as arbitrary.
The new feature that is introduced in this paper is that 
we also distinguish two kinds of equations:
strong equations and \emph{semi-equations},  
respectively denoted $\eqst$ and $\eqwe$,
so that equations also are \emph{decorated}.
Strong equations can be seen, essentially, as equalities between computations,
while semi-equations are much weaker,
and can be seen as a kind of approximation relation.
Moreover, as suggested by the symbols $\eqst$ and $\eqwe$,
the strong equations form an equivalence relation,
while the semi-equations form a preorder relation.
Then, we define the \emph{semi-product} of two functions
when at least one is pure, 
by a characteristic property that is a decorated version 
of the characteristic property of the usual product.
Since all identities are values, we get 
the semi-product of any function with an identity, 
that is used for building sequential products of functions.

Cartesian effect categories give rise to Arrows, in the sense of \cite{Hughes00},
and they provide a deduction system:
it is possible to decorate many proofs on cartesian categories 
in order to get proofs on cartesian effect categories.

\medskip

As for terminology, our \emph{graphs} are directed multi-graphs, 
made of \emph{points} (or vertices, or objects)
and \emph{functions} (or edges, arrows, morphisms).
We use \emph{weak} categories rather than categories, i.e., 
we use a congruence $\equiv$ rather than the equality,
however this ``syntactic'' choice is not fundamental here.
As for notations, 
we often omit the subscripts in the diagrams and in the proofs.

\medskip

Cartesian weak categories are reminded in section~\ref{sec:weak},
then cartesian effect categories are defined in section~\ref{sec:effect};
they are compared with Arrows in section~\ref{sec:related}, 
and examples are presented in section~\ref{sec:exam}.
In appendix~\ref{app:proof} are given the proofs 
of some properties of cartesian weak categories, that are well-known,
followed by their decorated versions, 
that yield proofs of properties of cartesian effect categories.


\section{Cartesian weak categories} 
\label{sec:weak}

Weak categories are reminded in this section, 
with their notion of product.
Except for the minor fact that equality is weakened as a congruence,
all this section is very well known.
Some detailed proofs are given in appendix~\ref{app:proof}, 
with their decorated versions. 


\subsection{Weak categories} 
\label{subsec:weak-cat}

A weak category is like a category,
except that the equations (for unitarity and associativity)
hold only ``up to congruence''. 
\begin{defi}
\label{defi:weak-cat}
A \emph{weak category} is a graph where: 
\begin{itemize}
\item for each point $X$ 
there is a loop $\id_X:X\to X$
called the \emph{identity} of $X$, 
\item for each consecutive functions $f:X\to Y$, $g:Y\to Z$, 
there is a function $g\circ f:X\to Z$
called the \emph{composition} of $f$ and $g$, 
\item and there is a relation $\equiv$ between 
parallel functions (each $f_1\equiv f_2$ is called an \emph{equation}),
such that:
\begin{itemize}
\item $\equiv$ is a \emph{congruence}, i.e., 
it is an equivalence relation and 
for each $f:X\to Y$, $g_1,g_2:Y\to Z$, $h:Z\to W$,
if $g_1\equiv g_2$ then $g_1\circ f\equiv g_2\circ f$ (\emph{substitution})
and $h\circ g_1 \equiv h\circ g_2$ (\emph{replacement}),
\item for each $f:X\to Y$, 
the \emph{unitarity equations} hold: 
$f\circ\id_X\equiv f$ and $\id_Y\circ f\equiv f$, 
\item and for each $f:X\to Y$, $g:Y\to Z$, $h:Z\to W$,
the \emph{associativity equation} holds: 
$h\circ (g \circ f) \equiv (h\circ g) \circ f$.
\end{itemize}
\end{itemize}
\end{defi}
\noindent So, a weak category is a special kind of a bicategory,
and a category is a weak category where the congruence is the equality. 


\subsection{Products}
\label{subsec:weak-prod}

In a weak category, a \emph{weak product}, or simply a \emph{product}, 
is defined as a product ``up to congruence''.
We focus on nullary products (i.e., terminal points)
and binary products; it is well-know that products of any arity 
can be recovered from those. 
\begin{defi}
\label{defi:weak-term}
A \emph{(weak) terminal point} is a point $U$ (for ``\emph{Unit}'')
such that for every point $X$ there is a function $\np_X:X\to U$, 
unique up to congruence.
\end{defi}
\begin{defi} 
\label{defi:weak-prod}
A \emph{binary cone} is made of two functions
with the same source $\cone{Y_1}{f_1}{X}{f_2}{Y_2}$.
A \emph{binary (weak) product} is a binary cone 
$\cone{Y_1}{q_1}{Y_1\times Y_2}{q_2}{Y_2}$
such that for every binary cone with the same base 
$\cone{Y_1}{f_1}{X}{f_2}{Y_2}$
there is a function 
$\tuple{f_1,f_2}:X\to Y_1\times Y_2$, 
called the \emph{pair} of $f_1$ and $f_2$, 
unique up to congruence, such that:
$$ q_1\circ \tuple{f_1,f_2}\equiv f_1  \;\mbox{ and }\;
 q_2\circ \tuple{f_1,f_2}\equiv f_2 \;.$$
\end{defi}
As usual, all terminal points are isomorphic, and 
the fact of using $U$ for denoting a terminal point 
corresponds to the choice of one terminal point. 
Similarly, all products on a given base are  isomorphic (in a suitable sense),
and the notations correspond to the choice of one product for each base.
\begin{defi}
\label{defi:weak-ccat}
A \emph{cartesian weak category} 
is a weak category with a chosen terminal point 
and chosen binary products. 
\end{defi}


\subsection{Products of functions}
\label{subsec:weak-prodfn}

\begin{defi}
\label{defi:weak-arr-prod}
In a cartesian weak category,
the \emph{(weak) binary product} of two functions 
$f_1:X_1\to Y_1$ and $f_2:X_2\to Y_2$ is the function:
$$ f_1\times f_2 = \tuple{f_1\circ p_1,f_2\circ p_2}:
X_1\times X_2 \to Y_1\times Y_2 \;.$$
\end{defi}
\noindent So, the binary product of functions is characterized, up to congruence, 
by the equations:
$$ q_1 \circ (f_1\times f_2) \equiv  f_1 \circ p_1 \;\mbox{ and }\; 
  q_2 \circ (f_1\times f_2) \equiv  f_2 \circ p_2 \;.$$
The defining equations of a pair and a product can be illustrated as follows:
  $$ 
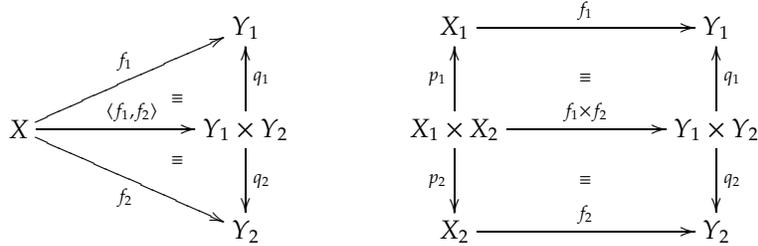
\xymatrix@C=5pc{ 
  & Y_1  \\
  X \ar[ru]^{f_1}  \ar[rd]_{f_2} \ar[r]^{\tuple{f_1,f_2}} &
    Y_1 \times Y_2 \ar[u]_{q_1} \ar[d]^{q_2}  
    \ar@{}[ld]|(.3){\equiv} \ar@{}[lu]|(.3){\equiv} \\
  & Y_2  \\
  }
  \qquad\qquad
  \xymatrix@C=5pc{ 
  X_1 \ar[r]^{f_1} & Y_1  \\
  X_1\times X_2 \ar[u]^{p_1} \ar[d]_{p_2} \ar[r]^{f_1\times f_2} 
    \ar@{}[rd]|{\equiv} \ar@{}[ru]|{\equiv } & 
    Y_1 \times Y_2 \ar[u]_{q_1} \ar[d]^{q_2}  \\
  X_2 \ar[r]^{f_2} & Y_2 \\
  } $$
So, the products are defined from the pairs
(note that we use the same symbols $f_1,f_2$ for 
the general case $f_i:X_i\to Y_i$ and for the special case $f_i:X\to Y_i$).
The other way round, the pairs can be recovered from the products
and the \emph{diagonals}, i.e., the pairs $\tuple{\id,\id}$;
indeed, it is easy to prove that for each cone $\cone{X_1}{f_1}{X}{f_2}{X_2}$
  $$ \tuple{f_1,f_2} \equiv (f_1\times f_2) \circ \tuple{\id_X,\id_X}\;.$$

In the following, we consider products
$\cone{X_1}{p_1}{X_1\times X_2}{p_2}{X_2}$,
$\cone{Y_1}{q_1}{Y_1\times Y_2}{q_2}{Y_2}$
and $\cone{Z_1}{r_1}{Z_1\times Z_2}{r_2}{Z_2}$.
\begin{prop}[congruence]
\label{prop:weak-equiv}
For each $f_1\equiv f'_1:X_1\to Y_1$ and $f_2\equiv f'_2:X_2\to Y_2$
\begin{enumerate}
\item if $X_1=X_2$
  $$ \tuple{f_1,f_2} \equiv \tuple{f'_1,f'_2} \;,$$
\item in all cases
  $$ f_1 \times f_2 \equiv f'_1 \times f'_2 \;.$$
\end{enumerate}
\end{prop}

\begin{prop}[composition]
\label{prop:weak-comp}
For each $f_1:X_1\to Y_1$, $f_2:X_2\to Y_2$, 
$g_1:Y_1\to Z_1$, $g_2:Y_2\to Z_2$
\begin{enumerate}
\item if $X_1=X_2$ and $Y_1=Y_2$ and $f_1=f_2(=f)$
  $$ \tuple{g_1,g_2} \circ f \equiv \tuple{g_1\circ f,g_2\circ f} \;,$$
\item if $X_1=X_2$
  $$ (g_1 \times g_2) \circ \tuple{f_1,f_2} \equiv 
      \tuple{g_1\circ f_1,g_2\circ f_2} \;,$$ 
\item in all cases
  $$ (g_1 \times g_2) \circ (f_1 \times f_2) \equiv 
   (g_1\circ f_1) \times (g_2\circ f_2) \;.$$ 
\end{enumerate}
\end{prop}

Let us consider the products
$\cone{X_1}{p_1}{X_1\times X_2}{p_2}{X_2}$  
and $\cone{X_2}{p'_2}{X_2\times X_1}{p'_1}{X_1}$.
The \emph{swap} function is the isomorphism:
  $$ \gamma_{(X_1,X_2)} = \tuple{p'_1,p'_2}_{p_1,p_2} = \tuple{p'_1,p'_2}:
   X_2\times X_1\to X_1\times X_2\;,$$
characterized by:
  $$ p_1\circ \gamma_{(X_1,X_2)}\equiv p'_1
  \;\mbox{ and }\; p_2\circ \gamma_{(X_1,X_2)}\equiv p'_2 \;.$$
\begin{prop}[swap]
\label{prop:weak-swap}
For each $f_1:X_1\to Y_1$ and $f_2:X_2\to Y_2$,
let $\gamma_Y=\gamma_{(Y_1,Y_2)}$ and  $\gamma_X=\gamma_{(X_1,X_2)}$,
then:
\begin{enumerate}
\item if $X_1=X_2$
  $$ \gamma_Y \circ \tuple{f_2,f_1} \equiv \tuple{f_1,f_2} \;,$$
\item in all cases
  $$ \gamma_Y \circ (f_2\times f_1) \circ \gamma_X^{-1} \equiv 
     f_1 \times f_2 \;.$$
\end{enumerate}
\end{prop}

Let us consider the products  
$\cone{X_1}{p_1}{X_1\times X_2}{p_2}{X_2}$,
$\cone{X_1\times X_2}{p_{1,2}}{(X_1\times X_2)\times X_3}{p_3}{X_3}$,
$\cone{X_2}{p'_2}{X_2\times X_3}{p'_3}{X_3}$ and 
$\cone{X_1}{p'_1}{X_1\times (X_2\times X_3)}{p'_{2,3}}{X_2\times X_3}$.
The \emph{associativity} function is the isomorphism:
  $$ \alpha_{(X_1,X_2,X_3)} = 
  \tuple{\tuple{p'_1,p'_2\circ p'_{2,3}}_{p_1,p_2},
    p'_3\circ p'_{2,3}}_{p_{1,2},p_3}:
   X_1\times (X_2\times X_3) \to (X_1\times X_2)\times X_3\;,$$
characterized by:
  $$ p_1\circ p_{1,2}\circ \alpha_{(X_1,X_2,X_3)} \equiv p'_1
  \,,\; 
  p_2\circ p_{1,2}\circ \alpha_{(X_1,X_2,X_3)} \equiv p'_2\circ p'_{2,3}
  \;\mbox{ and }\; 
  p_3\circ \alpha_{(X_1,X_2,X_3)} \equiv p'_3\circ p'_{2,3} \;.$$
\begin{prop}[associativity]
\label{prop:weak-assoc}
For each $f_1:X_1\to Y_1$, $f_2:X_2\to Y_2$ and $f_3:X_3\to Y_3$,
let $\alpha_Y=\alpha_{(Y_1,Y_2,Y_3)}$ 
and $\alpha_X=\alpha_{(X_1,X_2,X_3)}$,
then: 
\begin{enumerate}
\item if $X_1=X_2=X_3$
  $$ \alpha_Y \circ \tuple{f_1,\tuple{f_2,f_3}} 
  \equiv \tuple{\tuple{f_1,f_2},f_3} \;,$$
\item in all cases
  $$ \alpha_Y \circ (f_1\times (f_2\times f_3)) \equiv 
     ((f_1 \times f_2)\times f_3) \circ \alpha_X \;.$$
\end{enumerate}
\end{prop}

In the definition of the binary product $f_1\times f_2$,
both $f_1$ and $f_2$ play symmetric r\^oles.
This symmetry can be broken:
``first $f_1$ then $f_2$'' corresponds to 
$(\id_{Y_1}\times f_2) \circ (f_1\times\id_{X_2})$, 
using the intermediate product $Y_1\times X_2$,
while ``first $f_2$ then $f_1$'' corresponds to 
$(f_1\times\id_{Y_2}) \circ (\id_{X_1}\times f_2)$,
using the intermediate product $X_1\times Y_2$.
These are called the \emph{(left and right) sequential products}
of $f_1$ and $f_2$. 
The three versions of the binary product of functions coincide, 
up to congruence; this is a kind of \emph{parallelism} property,
meaning that both $f_1$ and $f_2$ can be computed either simultaneously, 
or one after the other, in any order:
\begin{prop}[parallelism]
\label{prop:weak-seq}
For each $f_1:X_1\to Y_1$ and $f_2:X_2\to Y_2$
$$ f_1\times f_2 \equiv (\id_{Y_1}\times f_2) \circ (f_1\times\id_{X_2})
\equiv (f_1\times\id_{X_2}) \circ (\id_{Y_1}\times f_2) \;.$$
\end{prop}


\section{Cartesian effect categories}
\label{sec:effect}

Sections~\ref{subsec:effect-cat} to~\ref{subsec:effect-semifn} 
form a decorated version of section~\ref{sec:weak}.
Roughly speaking, a kind of structure is \emph{decorated} 
when there is some classification of its ingredients.
Here, the classification involves two kinds of functions 
and two kinds of equations.
Effect categories are defined in section~\ref{subsec:effect-cat}
as decorated weak categories.
In section~\ref{subsec:effect-semi}, semi-products 
are defined as decorated weak products,
then cartesian effect category 
as decorated cartesian weak categories.
Decorated propositions are stated here, 
and the corresponding decorated proofs are given 
in appendix~\ref{app:proof}.
Then, in sections~\ref{subsec:effect-seq} and~\ref{subsec:effect-proj}, 
the sequential product of functions is defined by composing semi-products,
and some of its properties are derived.


\subsection{Effect categories} 
\label{subsec:effect-cat} 

A \emph{(weak) subcategory} $\bV$ of a weak category $\bC$ 
is a subcategory of $\bC$ such that each equation of $\bV$
is an equation of $\bC$. It is a \emph{wide (weak) subcategory} 
when $\bV$ and $\bC$ have the same points, 
and each equation of $\bC$ between functions in $\bV$
is an equation in $\bV$. 
Then only one symbol $\eqst$ can be used, for both 
$\bV$ and $\bC$.
\begin{defi}
\label{defi:effect-cat}
Let $\bV$ be a weak category.
An \emph{effect category extending} $\bV$ 
is a weak category $\bC$,
such that $\bV$ is a wide subcategory of $\bC$, 
together with a relation $\eqwe$ between parallel functions in $\bC$ such that:
\begin{itemize}
\item the relation $\eqwe$ is weaker than $\eqst$
for $f_1,f_2$ in $\bC$, $\,f_1\eqst f_2 \Rightarrow f_1\eqwe f_2$;
\item $\eqwe$ is transitive;
\item $\eqwe$ and $\eqst$ coincide on $\bV$
for $v_1,v_2$ in $\bV$, $\,v_1\eqst v_2 \iff v_1\eqwe v_2$;
\item $\eqwe$ satisfies the substitution property: 
\\ if $f:X\to Y$ and $g_1\eqwe g_2:Y\to Z$ then 
$g_1\circ f\eqwe g_2\circ f:X\to Z$;
\item $\eqwe$ satisfies the replacement property 
\emph{with respect to $\bV$}: 
\\ if $g_1\eqwe g_2:Y\to Z$ and $v:Z\to W$ in $\bV$
then $v\circ g_1 \eqwe v\circ g_2:Y\to W$.
\end{itemize}
\end{defi}
\noindent The first property implies that $\eqwe$ is reflexive,
and when $\eqst$ is the equality it means precisely 
that $\eqwe$ is reflexive. 
Since $\eqwe$ is transitive and weaker than $\eqst$, 
if either $f_1\eqst f_2\eqwe f_3$ or $f_1\eqwe f_2\eqst f_3$,
then $f_1\eqwe f_3$; this is called the \emph{compatibility}
of $\eqwe$ with $\eqst$.
An effect category is \emph{strict} when $\eqst$ is the equality. 
In this paper, 
there is no major difference between effect categories
and strict effect categories.

A \emph{pure function} is a function in $\bV$.
The symbol $\sto$ is used for pure functions, 
and $\to$ for all functions. 
It follows from definition~\ref{defi:effect-cat} that 
all the identities of $\bC$ are pure, 
the composition of pure functions is pure, 
and more precisely a composition of functions is pure 
if and only if all the composing functions are pure. 
It should be noted that there can be equations $f\eqst v$ 
between a non-pure function and a pure one;
then the function $f$ is proved effect-free, without being pure.
This ``syntactic'' choice could be argued; note that 
this situation disappears when the congruence $\eqst$ is the equality.
The relation $\eqwe$ is called the \emph{semi-congruence}
of the effect category, and each $f_1\eqwe f_2$ is called a \emph{semi-equation}.
The semi-congruence generally is not a congruence, for two reasons:
it may not be symmetric, and it may not 
satisfy the replacement property for all functions.

Examples of strict effect categories are given in section~\ref{sec:exam}.
For dealing with partiality in section~\ref{subsec:exam-partial},
the semi-congruence $\eqwe$ coincides with the usual ordering of partial functions,
it is not symmetric but it satisfies the replacement property for all partial functions.
On the other hand, in section~\ref{subsec:exam-state}, the semi-congruence $\eqwe$ 
means that two functions in an imperative language have the same   
result but may act differently on the state,
it is an equivalence relation that does not satisfy the replacement property
for non-pure functions. 

Clearly, if the decorations are forgotten, i.e.,
if both the distinction between pure functions and arbitrary functions 
and the distinction between the congruence and the semi-congruence 
are forgotten, then an effect category is just a weak category.

A cartesian effect category, as defined below, 
is an effect category where $\bV$ is cartesian 
and where this cartesian structure on $\bV$ 
has some kind of generalization to $\bC$,
that does \emph{not}, in general, turn $\bC$ into 
a cartesian weak category.


\subsection{Semi-products}
\label{subsec:effect-semi}

Now, let us assume that $\bC$ is an effect category
extending $\bV$, and that $\bV$ is cartesian.
We define nullary and binary \emph{semi-products} in $\bC$,
for building pairs of functions when at least one of them is pure.
\begin{defi}
\label{defi:effect-term}
A \emph{semi-terminal point} in $\bC$ is a terminal point $U$ in $\bV$ 
such that every function $g:X\to U$ satisfies $g\eqwe \np_X$.
\end{defi}
\begin{defi}
\label{defi:effect-prod}
A \emph{binary semi-product} in $\bC$ is a binary product 
$\scones{Y_1}{q_1}{Y_1\times Y_2}{q_2}{Y_2}$ in $\bV$ 
such that:
\begin{itemize}
\item for every binary cone with the same base 
$\cones{Y_1}{f_1}{X}{v_2}{Y_2}$ and with $v_2$ pure, 
there is a function
$\tuple{f_1,v_2}_{q_1,q_2}=\tuple{f_1,v_2}:X\to Y_1\times Y_2$, 
unique up to $\eqst$, such that 
$$ q_1\circ \tuple{f_1,v_2} \eqst f_1 \;\mbox{ and }\;  
  q_2\circ \tuple{f_1,v_2} \eqwe v_2 \;,$$
\item and for every binary cone with the same base 
$\scone{Y_1}{v_1}{X}{f_2}{Y_2}$ and with $v_1$ pure, 
there is a function
$\tuple{v_1,f_2}_{q_1,q_2}=\tuple{v_1,f_2}:X\to Y_1\times Y_2$, 
unique up to $\eqst$, such that 
$$ q_1\circ \tuple{v_1,f_2} \eqwe v_1 \;\mbox{ and }\;  
  q_2\circ \tuple{v_1,f_2} \eqst f_2 \;.$$
\end{itemize}
\end{defi}
\noindent The defining (semi-)equations of a binary semi-product can be illustrated 
as follows:
  $$ 
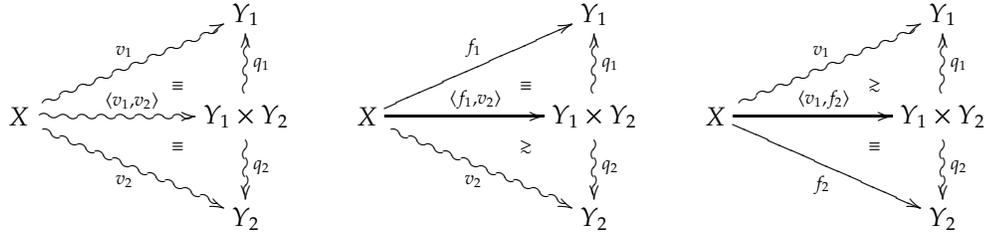
\xymatrix@C=5pc{ 
  & Y_1  \\
  X \ar@{~>}[ru]^{v_1} \ar@{~>}[rd]_{v_2} \ar@{~>}[r]^{\tuple{v_1,v_2}} &
    Y_1 \times Y_2 \ar@{~>}[u]_{q_1} \ar@{~>}[d]^{q_2}  
    \ar@{}[ld]|(.3){\eqst } \ar@{}[lu]|(.3){\eqst }\\
  & Y_2  \\
  }
  \qquad 
  \xymatrix@C=5pc{ 
  & Y_1  \\
  X \ar[ru]^{f_1}  \ar@{~>}[rd]_{v_2} \ar[r]^{\tuple{f_1,v_2}\;\;} &
    Y_1 \times Y_2 \ar@{~>}[u]_{q_1} \ar@{~>}[d]^{q_2}  
    \ar@{}[ld]|(.3){\weeq} \ar@{}[lu]|(.3){\eqst}\\
  & Y_2  \\
  } 
  \qquad
  \xymatrix@C=5pc{ 
  & Y_1  \\
  X \ar@{~>}[ru]^{v_1}  \ar[rd]_{f_2} \ar[r]^{\tuple{v_1,f_2}\;\;} &
    Y_1 \times Y_2 \ar@{~>}[u]_{q_1} \ar@{~>}[d]^{q_2}  
    \ar@{}[ld]|(.3){\eqst} \ar@{}[lu]|(.3){\weeq}\\
  & Y_2  \\
  }   
$$
Clearly, if the decorations are forgotten, 
then semi-products are just products.

The notation is not ambiguous. Indeed, if 
$\scones{Y_1}{v_1}{X}{v_2}{Y_2}$ is a binary cone in $\bV$, 
then the three definitions of the pair $\tuple{v_1,v_2}$
above coincide, up to congruence:
let $t$ denote any one of the three pairs, then 
$t$ is characterized, up to congruence, 
by $q_1\circ t \eqst v_1$ and $q_2\circ t \eqst v_2$,
because $\eqst$ and $\eqwe$ coincide on pure functions.

\begin{defi}
\label{defi:effect-ccat}
A \emph{cartesian effect category extending} a cartesian weak 
category $\bV$ 
is an effect category extending $\bV$ 
such that each terminal point of $\bV$ is a semi-terminal point of $\bC$ 
and each binary product of $\bV$ is a binary semi-product of $\bC$.
\end{defi}


\subsection{Semi-products of functions}
\label{subsec:effect-semifn}

\begin{defi}
\label{defi:effect-arr-prod}
In a cartesian effect category,
the \emph{binary semi-product} $f_1\times v_2$ of 
a function $f_1:X_1\to Y_1$ and a pure function $v_2:X_2\sto Y_2$ is
the function:
$$f_1\times v_2=\tuple{f_1\circ p_1,v_2\circ p_2} : 
  X_1\times X_2 \to Y_1\times Y_2$$
\end{defi}
\noindent It follows that $f_1\times v_2$ is characterized, up to $\eqst$, 
by: 
$$ q_1 \circ (f_1\times v_2) \eqst  f_1 \circ p_1 \;\mbox{ and }\; 
  q_2 \circ (f_1\times v_2) \eqwe  v_2 \circ p_2 $$
$$ \xymatrix@C=5pc{ 
  X_1 \ar[r]^{f_1} & Y_1  \\
  X_1\times X_2 \ar@{~>}[u]^{p_1} \ar@{~>}[d]_{p_2} \ar[r]^{f_1\times v_2} 
    \ar@{}[rd]|{\weeq} \ar@{}[ru]|{\eqst } & 
    Y_1 \times Y_2 \ar@{~>}[u]_{q_1} \ar@{~>}[d]^{q_2}  \\
  X_2 \ar@{~>}[r]^{v_2} & Y_2 \\
  } $$
The \emph{binary semi-product} 
$v_1\times f_2:X_1\times X_2\to Y_1\times Y_2$ of 
a pure function $v_1:X_1\sto Y_1$ and a function $f_2:X_2\to Y_2$ 
is defined in the symmetric way, and it is characterized, up to $\eqst$, 
by the symmetric property.

The notation is not ambiguous, because so is the notation for pairs; 
if $v_1$ and $v_2$ are pure functions, 
then the three definitions of $v_1\times v_2$ coincide, up to congruence.

Propositions about products in cartesian weak categories 
are called \emph{basic} propositions.
It happens that each basic proposition in section~\ref{sec:weak} 
has a \emph{decorated} version, 
about semi-products of the form $f_1\times v_2$ in cartesian effect categories,
that is stated below. The symmetric decorated version also holds, 
for semi-products of the form $v_1\times f_2$.
Each function in the basic proposition
is replaced either by a function or by a pure function,
and each equation 
is replaced either by an equation ($\eqst$) or by a semi-equation
($\eqwe$ or $\weeq$).

In addition, in appendix~\ref{app:proof},
the proofs of the decorated propositions  
are \emph{decorated} versions of the \emph{basic} proofs.
It happens that no semi-equation appears in the 
decorated propositions below, but they are used in the proofs.
Indeed, a major ingredient in the basic proofs is that 
a function $\tuple{f_1,f_2}$ or $f_1\times f_2$ 
is characterized, up to $\equiv$,
by its projections, both up to $\equiv$.
The decorated version of this property is that 
a function $\tuple{f_1,f_2}$ or $f_1\times f_2$,
where $f_1$ or $f_2$ is pure, 
is characterized, up to $\equiv$,
by its projections, one up to $\equiv$ \emph{and the other one 
up to $\eqwe$}. 
It should be noted that even when some decorated version of a 
basic proposition is valid, 
usually not all the basic proofs can be decorated.
In addition, when equations are decorated as semi-equations,
some care is required 
when the symmetry and replacement properties are used.

\begin{prop}[congruence]
\label{prop:effect-equiv}
For each congruent functions $f_1\eqst f'_1:X\to Y_1$
and pure functions $v_2\eqst v'_2:X\sto Y_2$
\begin{enumerate}
\item if $X_1=X_2$
  $$ \tuple{f_1,v_2} \eqst \tuple{f'_1,v'_2} \;.$$
\item in all cases
  $$ f_1 \times v_2 \eqst f'_1 \times v'_2 \;.$$
\end{enumerate}
\end{prop}

\begin{prop}[composition]
\label{prop:effect-comp}
For each functions $f_1:X_1\to Y_1$, $g_1:Y_1\to Z_1$
and pure functions $v_2:X_2\sto Y_2$, $w_2:Y_2\sto Z_2$
\begin{enumerate}
\item if $X_1=X_2$ and $Y_1=Y_2$ and $f_1=v_2(=v)$
  $$ \tuple{g_1,w_2} \circ f \eqst \tuple{g_1\circ v,w_2\circ v} \;,$$
\item if $X_1=X_2$
  $$ (g_1 \times w_2) \circ \tuple{f_1,v_2} \eqst 
      \tuple{g_1\circ f_1,w_2\circ v_2} \;,$$ 
\item in all cases
  $$ (g_1 \times w_2) \circ (f_1 \times v_2) \eqst 
   (g_1\circ f_1) \times (w_2\circ v_2) \;.$$ 
\end{enumerate}
\end{prop}

The swap and associativity functions are defined in the same
way as in section~\ref{sec:weak}; 
they are products of projections, so that they are pure functions.
It follows that the swap and associativity functions
are characterized by the same equations as in section~\ref{sec:weak}, 
and that they are still isomorphisms.

\begin{prop}[swap]
\label{prop:effect-swap}
For each function $f_1:X\to Y_1$ and pure function $v_2:X\sto Y_2$,
let $\gamma_Y=\gamma_{(Y_1,Y_2)}$ and  $\gamma_X=\gamma_{(X_1,X_2)}$,
then:
\begin{enumerate}
\item if $X_1=X_2$
  $$ \gamma_Y \circ \tuple{v_2,f_1} \eqst \tuple{f_1,v_2} \;,$$
\item in all cases
  $$ \gamma_Y \circ (v_2\times f_1) \circ \gamma_X^{-1} \eqst 
     f_1 \times v_2 \;.$$
\end{enumerate}
\end{prop}

\begin{prop}[associativity]
\label{prop:effect-assoc}
For each function $f_1:X_1\to Y_1$
and pure functions $v_2:X_2\sto Y_2$, $v_3:X_3\sto Y_3$,
let $\alpha_Y=\alpha_{(Y_1,Y_2,Y_3)}$ 
and $\alpha_X=\alpha_{(X_1,X_2,X_3)}$,
then: 
\begin{enumerate}
\item if $X_1=X_2=X_3$
  $$ \alpha_Y \circ \tuple{f_1,\tuple{v_2,v_3}} 
  \eqst \tuple{\tuple{f_1,v_2},v_3} \;,$$
\item in all cases:
  $$ \alpha_Y \circ (f_1\times (v_2\times v_3)) \eqst 
     ((f_1 \times v_2)\times v_3) \circ \alpha_X \;.$$
\end{enumerate}
\end{prop}

The sequential product of a function $f_1:X_1\to Y_1$ 
and a pure function $v_2:X_2\sto Y_2$ can be defined 
as in section~\ref{sec:weak}, 
using the intermediate
products $\scones{Y_1}{s_1}{Y_1\times X_2}{s_2}{X_2}$
and $\scones{X_1}{t_1}{X_1\times Y_2}{t_2}{Y_2}$. 
It does coincide with the semi-product of $f_1$ and $v_2$, 
up to congruence:
\begin{prop}[parallelism]
\label{prop:effect-seq}
For each function $f_1:X_1\to Y_1$ and pure function $v_2:X_2\sto Y_2$
$$ f_1\times v_2 \eqst (\id_{Y_1}\times v_2) \circ (f_1\times\id_{X_2})
\eqst (f_1\times\id_{X_2}) \circ (\id_{Y_1}\times v_2) \;.$$
\end{prop}


\subsection{Sequential products of functions}
\label{subsec:effect-seq}

It has been stated in proposition~\ref{prop:weak-seq}
that, in a cartesian weak category,
the binary product of functions coincide with both
sequential products, up to congruence:
$$ f_1\times f_2 \equiv (\id_{Y_1}\times f_2) \circ (f_1\times\id_{X_2})
\equiv (f_1\times\id_{X_2}) \circ (\id_{Y_1}\times f_2) \;.$$
In a cartesian effect category, 
when $f_1$ and $f_2$ are any functions, 
the product $f_1\times f_2$ is not defined.
But $(\id_{Y_1}\times f_2) \circ (f_1\times\id_{X_2})$
and $(f_1\times\id_{X_2}) \circ (\id_{Y_1}\times f_2)$
make sense, thanks to semi-products,
because identities are pure. 
They are called the sequential products of $f_1$ and $f_2$,
and they do not coincide up to congruence, in general:
parallelism is not satisfied.
\begin{defi}
\label{defi:effect-seq-prod}
The \emph{left binary sequential product} of two functions
$f_1:X_1\to Y_1$ and $f_2:X_2\to Y_2$ is the function:
$$ f_1\ltimes f_2 = 
  (\id_{Y_1}\times f_2) \circ (f_1\times\id_{X_2}) : X_1\times X_2\to Y_1\times Y_2 \;.$$
\end{defi}
\noindent So, the left binary sequential product is obtained from:
$$ \xymatrix@C=5pc{ 
  X_1 \ar[r]^{f_1} & Y_1 \ar@{~>}[r]^{\id} & Y_1  \\
  X_1\times X_2 \ar@{~>}[u]^{p_1} \ar@{~>}[d]_{p_2} 
    \ar[r]^{f_1\times \id} 
    \ar@{}[rd]|{\weeq} \ar@{}[ru]|{\eqst} & 
    Y_1 \times X_2 \ar@{~>}[u]^{s_1} \ar@{~>}[d]_{s_2}  
    \ar[r]^{\id\times f_2} 
    \ar@{}[rd]|{\eqst} \ar@{}[ru]|{\weeq} & 
    Y_1 \times Y_2 \ar@{~>}[u]^{q_1} \ar@{~>}[d]_{q_2}  \\
  X_2 \ar@{~>}[r]^{\id} & X_2 \ar[r]^{f_2} & Y_2 \\
  } $$
The left sequential product extends the semi-product: 
\begin{prop}
\label{prop:effect-seq-lproduct}
For each function $f_1$ and pure function $v_2$, 
$ f_1\ltimes v_2 \eqst f_1\times v_2 $.
\end{prop}
\begin{proof}
>From proposition~\ref{prop:effect-comp}, 
$f_1\ltimes v_2 = (\id\times v_2) \circ (f_1\times\id)
\eqst  (\id\circ f_1) \times (v_2\circ\id)
\eqst f_1 \times v_2$.
\end{proof}
\noindent Note that the diagonal $\tuple{\id_X,\id_X}$ is a pair of pure functions. 
So, by analogy with the property  
$\tuple{f_1,f_2} \equiv (f_1\times f_2) \circ \tuple{\id_X,\id_X}$ 
in weak categories:
\begin{defi}
\label{defi:effect-seq-pair} 
The \emph{left sequential pair} of two functions 
$f_1:X\to Y_1$ and $f_2:X\to Y_2$ is: 
  $$ \tuple{f_1,f_2}_l = (f_1\ltimes f_2) \circ \tuple{\id_X,\id_X}\;.$$
\end{defi} 
\noindent The left sequential pairs do not satisfy the usual equations 
for pairs, as in definition~\ref{defi:weak-prod}.
However, they satisfy some weaker properties, 
as stated in corollary~\ref{cor:seq-pair}.

The \emph{right binary sequential product} of $f_1$ and $f_2$
is defined in the symmetric way; 
it is the function:
$$ f_1\rtimes f_2 = 
  (f_1\times\id_{Y_2}) \circ (\id_{X_1}\times f_2) : 
  X_1\times X_2\to Y_1\times Y_2 \;.$$
It does also extend the product of a pure function and a function:
for each pure function $v_1$, $ v_1\rtimes f_2 \eqst v_1\times f_2 $.
The  \emph{right sequential pair} of $f_1:X\to Y_1$ and $f_2:X\to Y_2$ is: 
  $$ \tuple{f_1,f_2}_r = (f_1\rtimes f_2) \circ \tuple{\id_X,\id_X}\;.$$

Here are some properties of the sequential products 
that are easily deduced from the properties 
of semi-products in~\ref{subsec:effect-semi}.
The symmetric properties also hold.

\begin{prop}[congruence]
\label{prop:seq-equiv}
For each congruent functions $f_1\eqst f'_1:X_1\to Y_1$
and $f_2\eqst f'_2:X_2\to Y_2$
  $$f_1 \ltimes f_2 \eqst f'_1 \ltimes f'_2\;.$$
\end{prop}
\begin{proof}
Clear, from~\ref{prop:effect-equiv}.
\end{proof}

\begin{prop}[composition]
\label{prop:seq-comp}
For each functions $f_1:X_1\to Y_1$, $g_1:Y_1\to Z_1$, $g_2:Y_2\to Z_2$
and pure function $v_2:X_2\sto Y_2$
  $$  (g_1 \ltimes g_2) \circ (f_1 \times v_2) \eqst 
   (g_1\circ f_1) \ltimes (g_2\circ v_2) \;.$$
\end{prop}
$$  \xymatrix@C=5pc{ 
  X_1 \ar[r]^{f_1} & Y_1 \ar[r]^{g_1} & Z_1 \ar@{~>}[r]^{\id} & Z_1  \\
  X_1 \times X_2 \ar@{~>}[u] \ar@{~>}[d]   \ar[r]^{f_1\times v_2} 
  \ar@{}[rd]|{\eqst} \ar@{}[ru]|{\eqst} &
  Y_1\times Y_2 \ar@{~>}[u] \ar@{~>}[d] 
    \ar[r]^{g_1\times \id} 
    \ar@{}[rd]|{\weeq} \ar@{}[ru]|{\eqst} & 
    Z_1 \times Y_2 \ar@{~>}[u] \ar@{~>}[d]  
    \ar[r]^{\id\times g_2} 
    \ar@{}[rd]|{\eqst} \ar@{}[ru]|{\weeq} & 
    Z_1 \times Z_2 \ar@{~>}[u] \ar@{~>}[d]  \\
  X_2 \ar@{~>}[r]^{v_2} & Y_2 \ar@{~>}[r]^{\id} & Y_2 \ar[r]^{g_2} & Z_2 \\
  } $$
\begin{proof}
>From several applications of proposition~\ref{prop:effect-comp}
and its symmetric version:
\\ $(\id\times g_2) \circ (g_1\times\id) \circ (f_1 \times v_2)
\eqst (\id\times g_2) \circ  ((g_1\circ f_1) \times v_2)
\eqst (\id\times g_2) \circ (\id\times v_2) 
\circ ((g_1\circ f_1)\times\id)
\eqst (\id\times (g_2\circ v_2)) \circ ((g_1\circ f_1)\times\id)$.
\end{proof}

\begin{prop}[swap]
\label{prop:seq-swap}
For each functions $f_1:X_1\to Y_1$ and $f_2:X_2\to Y_2$,
the left and right sequential products are related by swaps:
  $$ \gamma_Y \circ (f_2\rtimes f_1) \circ \gamma_X^{-1} \eqst 
     f_1 \ltimes f_2 \;.$$
\end{prop}
\begin{proof}
>From proposition~\ref{prop:effect-swap}
and its symmetric version:
\\ $ \gamma \circ (\id\times f_2) \circ (f_1\times\id)
\eqst  (f_2\times\id ) \circ \gamma \circ (f_1\times\id)
\eqst  (f_2\times\id ) \circ (\id\times f_1) \circ \gamma$.
\end{proof}

\begin{prop}[associativity]
\label{prop:seq-assoc}
For each functions $f_1:X_1\to Y_1$, $f_2:X_2\to Y_2$ and $f_3:X_3\to Y_3$,
let $\alpha_Y=\alpha_{(Y_1,Y_2,Y_3)}$ and $\alpha_X=\alpha_{(X_1,X_2,X_3)}$,
then: :
  $$ \alpha_Y \circ (f_1\ltimes (f_2\ltimes f_3)) \eqst 
     ((f_1 \ltimes f_2)\ltimes f_3) \circ \alpha_X \;.$$
\end{prop}
\begin{proof}
>From proposition~\ref{prop:effect-assoc}.
\end{proof}


\subsection{Projections of sequential products} 
\label{subsec:effect-proj} 

Let us come back to a weak category, as in section~\ref{sec:weak}. 
The binary product of functions is characterized, up to congruence, 
by the equations:
$$ q_1 \circ (f_1\times f_2) \equiv  f_1 \circ p_1 \;\mbox{ and }\; 
  q_2 \circ (f_1\times f_2) \equiv  f_2 \circ p_2 \;,$$
so that for all 
constant functions $x_1:U\to X_1$ and $x_2:U\to X_2$
$$ q_1 \circ (f_1\times f_2) \circ \tuple{x_1,x_2} \equiv f_1 \circ x_1 
\;\mbox{ and }\; 
  q_2 \circ (f_1\times f_2) \circ \tuple{x_1,x_2} \equiv f_2 \circ x_2 \;.$$
In a cartesian effect category, it is proved 
in theorem~\ref{thm:seq-prod} that $f_1\ltimes f_2$, 
when applied to a pair of constant pure functions $\tuple{x_1,x_2}$, 
returns on the $Y_1$ side a function that is semi-congruent to $f_1(x_1)$,
and on the $Y_2$ side a function that is congruent to 
$f_2\circ x_2\circ \np \circ f_1\circ x_1$,
which means ``first $f_1(x_1)$, then forget the result, then $f_2(x_2)$''. 
More precise statements are given in propositions~\ref{prop:seq-val}
and~\ref{prop:seq-com}.
Proofs are presented in the same formalized way as in appendix~\ref{app:proof}.

As above, we consider the semi-terminal point $U$
and semi-products 
\\ $\scones{X_1}{p_1}{X_1\times X_2}{p_2}{X_2}$,
$\scones{Y_1}{q_1}{Y_1\times Y_2}{q_2}{Y_2}$
and $\scones{Y_1}{s_1}{Y_1\times X_2}{s_2}{X_2}$.
\begin{prop}
\label{prop:seq-val}
For each functions $f_1:X_1\to Y_1$ and $f_2:X_2\to Y_2$
$$q_1\circ (f_1\ltimes f_2) \eqwe f_1\circ p_1 : X_1\times X_2\to Y_1\;.$$
\end{prop}
$$ \xymatrix@C=5pc{ 
  X_1 \ar[r]^{f_1} & Y_1 \ar@{~>}[r]^{\id} & Y_1  \\
  X_1\times X_2 \ar@{~>}[u]^{p_1} 
    \ar[r]^{f_1\times \id} 
    \ar@{}[ru]|{\eqst} \ar@/_4ex/[rr]_{f_1\ltimes f_2}^{=} & 
    Y_1 \times X_2 \ar@{~>}[u]^{s_1}   
    \ar[r]^{\id\times f_2} 
    \ar@{}[ru]|{\weeq} & 
    Y_1 \times Y_2 \ar@{~>}[u]^{q_1} \\
  } $$
\begin{proof}\
\\ \begin{tabular}{rll}
$(a)$ & $q_1\circ (\id\times f_2) \eqwe s_1$ & \\
$(b)$ & $q_1\circ (f_1\ltimes f_2) \eqwe s_1\circ (f_1\times\id)$ & $(a)$, $\subst_{\eqwe}$ \\
$(c)$ & $s_1\circ (f_1\times\id) \eqst f_1\circ p_1$ & \\
$(d)$ & $q_1\circ (f_1\ltimes f_2) \eqwe f_1\circ p_1$ & $(b)$, $(c)$, $\comp$ \\
\end{tabular}\\ 
\end{proof}
\begin{lem}
\label{lem:seq-terminal}
For each function $f_1:X_1\to Y_1$ and pure function $x_2:U\sto X_2$
$$ \tuple{\id_{Y_1},x_2\circ\np_{Y_1}}\circ f_1 \eqst \tuple{f_1,x_2\circ\np_{X_1}}
  :  X_1 \to Y_1\times X_2 \;.$$
\end{lem}
\noindent Both handsides can be illustrated as follows:
$$ \xymatrix@C=5pc{ 
  & Y_1 \ar@{~>}[r]^{\id} & Y_1  \\
  X_1 \ar[r]^{f_1} & Y_1 \ar@{=}[u] \ar@{~>}[d]_{\np}  
    \ar@{~>}[r]^{\tuple{\id,x_2\circ\np}}
    \ar@{}[rd]|{\eqst} \ar@{}[ru]|{\eqst} & 
    Y_1 \times X_2 \ar@{~>}[u]^{s_1} \ar@{~>}[d]_{s_2}  \\
  & U \ar@{~>}[r]^{x_2} & X_2 \\
  } \qquad \xymatrix@C=5pc{ 
  X_1 \ar[r]^{f_1} & Y_1  \\
  X_1 \ar@{=}[u] \ar@{~>}[d]_{\np}  
    \ar[r]^{\tuple{f_1,x_2\circ\np}}
    \ar@{}[rd]|{\weeq} \ar@{}[ru]|{\eqst} & 
    Y_1 \times X_2 \ar@{~>}[u]^{s_1} \ar@{~>}[d]_{s_2}  \\
  U \ar@{~>}[r]^{x_2} & X_2 \\
  } $$
\begin{proof}\
\\ \begin{tabular}{rll}
$(a_1)$ & $s_1\circ \tuple{\id,x_2\circ\np} \eqst \id$ & \\
$(b_1)$ & $s_1\circ \tuple{\id,x_2\circ\np}\circ f_1 \eqst f_1$ & $(a_1)$, $\subst_{\eqst}$ \\
$(a_2)$ & $s_2\circ \tuple{\id,x_2\circ\np} \eqst x_2\circ\np$ & \\
$(b_2)$ & $s_2\circ \tuple{\id,x_2\circ\np}\circ f_1 \eqst x_2\circ\np\circ f_1$ & $(a_2)$, $\subst_{\eqst}$ \\
$(c_2)$ & $\np\circ f_1 \eqwe\np$ & semi-terminality of $U$ \\
$(d_2)$ & $x_2\circ\np\circ f_1\eqwe x_2\circ\np$ &  $(c_2)$, $\repl_{\eqwe}$ ($x_2$ is pure) \\
$(e_2)$ & $s_2\circ \tuple{\id,x_2\circ\np}\circ f_1 \eqwe  x_2\circ\np$ & $(b_2)$, $(d_2)$, $\tr_{\eqwe}$ \\
$(f)$ & $\tuple{\id,x_2\circ\np}\circ f_1 \eqst \tuple{f_1,x_2\circ\np}$ & $(b_1)$, $(e_2)$ \\
\end{tabular}\\ 
\end{proof}
\begin{prop}
\label{prop:seq-com}
For each functions $f_1:X_1\to Y_1$, $f_2:X_2\to Y_2$ 
and pure function $x_2:U\sto X_2$
$$ q_2 \circ (f_1\ltimes f_2) \circ \tuple{\id_{X_1},x_2\circ \np_{X_1}}
\eqst f_2 \circ x_2 \circ \np_{Y_1} \circ f_1 : X_1\to Y_2 \;.$$
\end{prop}
\noindent Both handsides can be illustrated as follows:
  $$ \xymatrix@C=3pc{ 
  X_1 \ar@{~>}[r]^{\id} & X_1 & &  \\
  X_1 \ar@{=}[u] \ar@{~>}[d]_{\np} 
    \ar@{~>}[r]^{\tuple{\id,x_2\circ\np}\;} &
    X_1 \times X_2  \ar@/^4ex/[rr]^{f_1\ltimes f_2}_{=} 
    \ar@{~>}[u]_{p_1} \ar@{~>}[d]^{p_2} \ar[r]^{f_1\times\id}
    \ar@{}[ld]|{\eqst} \ar@{}[lu]|{\eqst} \ar@{}[rd]|{\weeq} & 
    Y_1 \times X_2 \ar@{~>}[d]^{s_2} \ar[r]^{\id\times f_2} 
    \ar@{}[rd]|{\eqst } & 
    Y_1 \times Y_2 \ar@{~>}[d]^{q_2}  \\
  U \ar@{~>}[r]^{x_2} & X_2 \ar@{~>}[r]^{\id} & X_2 \ar[r]^{f_2} & Y_2  \\
  } \quad 
  \xymatrix@C=2pc{ 
  \mbox{ } \\
  X_1 \ar[r]^{f_1} & Y_1 \ar@{~>}[d]_{\np}  \\
  & U \ar@{~>}[r]^{x_2} & X_2 \ar[r]^{f_2} & Y_2  \\
  } $$
\begin{proof}\
\\ \begin{tabular}{rll}
$(a)$ & $q_2\circ (\id\times f_2) \eqst f_2\circ s_2$ \\
$(b)$ & $q_2\circ (f_1\ltimes f_2) \circ  \tuple{\id,x_2\circ \np}
  \eqst f_2\circ s_2\circ (f_1\times\id) \circ \tuple{\id,x_2\circ \np}$ & $(a)$, $\subst_{\eqst}$ \\
$(c)$ & $(f_1\times\id) \circ \tuple{\id,x_2\circ \np} 
  \eqst \tuple{f_1,x_2\circ\np}$ & prop.~\ref{prop:effect-comp} \\
$(d)$ & $\tuple{f_1,x_2\circ\np}
  \eqst \tuple{\id,x_2\circ\np}\circ f_1$ & lemma~\ref{lem:seq-terminal} \\
$(e)$ & $(f_1\times\id) \circ \tuple{\id,x_2\circ \np} 
  \eqst \tuple{\id,x_2\circ\np}\circ f_1$ & $(c)$, $(d)$, $\tr_{\eqst}$ \\ 
$(f)$ & $f_2\circ s_2\circ (f_1\times\id) \circ \tuple{\id,x_2\circ \np}
  \eqst f_2\circ s_2\circ \tuple{\id,x_2\circ\np}\circ f_1$ & $(e)$, $\repl_{\eqst}$ \\ 
$(g)$ & $q_2\circ (f_1\ltimes f_2) \circ  \tuple{\id,x_2\circ \np}
  \eqst f_2\circ s_2\circ \tuple{\id,x_2\circ\np}\circ f_1$ & $(b)$, $(f)$, $\tr_{\eqst}$ \\ 
$(h)$ & $p_2\circ \tuple{\id,x_2\circ\np}  \eqst x_2\circ \np$ & \\
$(i)$ & $f_2\circ s_2\circ \tuple{\id,x_2\circ\np}\circ f_1
  \eqst f_2\circ x_2\circ \np \circ f_1 $ & $(h)$, $\subst_{\eqst}$, $\repl_{\eqst}$ \\ 
$(j)$ & $q_2\circ (f_1\ltimes f_2) \circ  \tuple{\id,x_2\circ \np}
  \eqst f_2\circ x_2\circ \np \circ f_1$ & $(g)$, $(i)$,  $\tr_{\eqst}$ \\ 
\end{tabular}\\ 
\end{proof}
\begin{thm}
\label{thm:seq-prod}
For each functions $f_1:X_1\to Y_1$, $f_2:X_2\to Y_2$ 
and pure functions $x_1:U\sto X_1$ and $x_2:U\sto X_2$,
the function $(f_1\ltimes f_2) \circ \tuple{x_1,x_2}$ satisfies:
$$ q_1 \circ (f_1\ltimes f_2) \circ \tuple{x_1,x_2} 
  \eqwe f_1 \circ x_1 \;\mbox{ and }\;
   q_2 \circ (f_1\ltimes f_2) \circ \tuple{x_1,x_2} 
  \eqst f_2 \circ x_2 \circ \np_{Y_1} \circ f_1 \circ x_1 \;.$$
\end{thm}
$$ \xymatrix@C=3pc{ 
  U \ar@{=}[d] \ar@{~>}[rr]^{x_1} && X_1 \ar[rrr]^{f_1} &&& Y_1 \\
  U \ar@{=}[d] \ar@{~>}[rr]^{\tuple{x_1,x_2}} && 
    X_1\times X_2 \ar[rrr]^{f_1\ltimes f_2} 
    \ar@{}[ru]|{\weeq} \ar@{}[rd]|{\eqst} &&& 
    Y_1 \times Y_2 \ar@{~>}[u]_{q_1} \ar@{~>}[d]^{q_2} \\
  U \ar@{~>}[r]^{x_1} & X_1 \ar[r]^{f_1} & Y_1 \ar@{~>}[r]^{\np}
  & U \ar@{~>}[r]^{x_2} & X_2 \ar[r]^{f_2} & Y_2 \\
  } $$
\begin{proof}\
\\ \begin{tabular}{rll}
$(a_1)$ & $q_1\circ (f_1\ltimes f_2) \eqwe f_1\circ p_1$ & prop.~\ref{prop:seq-val} \\ 
$(b_1)$ & $q_1\circ (f_1\ltimes f_2) \circ \tuple{x_1,x_2} 
\eqwe f_1\circ p_1 \circ \tuple{x_1,x_2}$ & $(a_1)$, $\subst_{\eqwe}$ \\
$(c_1)$ & $p_1 \circ \tuple{x_1,x_2}\eqst x_1$ &  (on values) \\
$(d_1)$ & $f_1\circ p_1 \circ \tuple{x_1,x_2}\eqst f_1\circ x_1$ & $(c_1)$, $\repl_{\eqst}$ \\
$(e_1)$ & $q_1 \circ (f_1\ltimes f_2) \circ \tuple{x_1,x_2} 
  \eqwe f_1 \circ x_1$ & $(b_1)$, $(d_1)$, $\comp$ \\
$(a_2)$ & $\tuple{x_1,x_2} \eqst \tuple{\id_{X_1},x_2\circ \np_{X_1}} \circ x_1$ 
  & (on values) \\
$(b_2)$ & $q_2 \circ (f_1\ltimes f_2) \circ \tuple{x_1,x_2} \eqst 
  q_2 \circ (f_1\ltimes f_2) \circ \tuple{\id_{X_1},x_2\circ \np_{X_1}} \circ x_1$ 
  & $(a_2)$, $\repl_{\eqst}$ \\
$(c_2)$ & $q_2 \circ (f_1\ltimes f_2) \circ \tuple{\id_{X_1},x_2\circ \np_{X_1}}
\eqst f_2 \circ x_2 \circ \np_{Y_1} \circ f_1$ & prop.~\ref{prop:seq-com} \\ 
$(d_2)$ & $q_2 \circ (f_1\ltimes f_2) \circ \tuple{\id_{X_1},x_2\circ \np_{X_1}} 
\circ x_1
\eqst f_2 \circ x_2 \circ \np_{Y_1} \circ f_1\circ x_1$ & $(c_2)$, $\subst_{\eqst}$  \\
$(e_2)$ & $q_2 \circ (f_1\ltimes f_2) \circ \tuple{x_1,x_2} \eqst 
f_2 \circ x_2 \circ \np_{Y_1} \circ f_1\circ x_1$ & $(b_2)$,  $(d_2)$, $\tr_{\eqst}$ \\
\end{tabular}\\
\end{proof}

The corresponding properties of left sequential pairs easily follow.
\begin{cor}
\label{cor:seq-pair}
For each functions $f_1:X\to Y_1$, $f_2:X\to Y_2$ 
and pure function $x:U\sto X$
  $$ q_1\circ \tuple{f_1,f_2}_l \eqwe f_1 \,\mbox{, hence }\,
  q_1\circ \tuple{f_1,f_2}_l\circ x \eqwe f_1\circ x 
  \,\mbox{, and }\, q_2\circ \tuple{f_1,f_2}_l\circ x \eqst  
  f_2 \circ x \circ \np_{Y_1} \circ f_1 \circ x \;.$$
\end{cor}


\section{Effect categories and Arrows}
\label{sec:related}



Starting from \cite{Moggi91,Wadler93}, \emph{monads} are used in Haskell 
for dealing with computational effects. 
A \emph{Monad type} in Haskell is a unary type constructor 
that corresponds to a \emph{strong monad}, in the categorical sense. 
Monads have been generalized 
on the categorical side to \emph{Freyd categories} \cite{PowerRobinson97}
and on the functional programming side to \emph{Arrows} \cite{Hughes00}. 
A precise statement of the facts 
that Arrows generalize Monads and that Arrows are Freyd categories
can be found in \cite{HeunenJacobs06},
where each of the three notions is seen as 
a monoid in a relevant category.
Now we prove that cartesian effect categories determine Arrows.
In section~\ref{sec:exam} our approach is compared with 
the Monads approach, for two fundamental examples.
In this section, all effect categories are strict:
the congruence $\eqst$ is the equality.


\subsection{Arrows} 
\label{subsec:related-arrows}

According to \cite{Paterson01}, Arrows in Haskell are defined as follows. 

\begin{defi}
\label{defi:related-arr}
An \emph{Arrow} is a binary type constructor class $\tA$ of the form:
\begin{tabbing}
999 \= 999 \= 999 \kill
\> \texttt{class Arrow $\tA$ where} \\
\> \> $\arr:: (X\to Y)\to \tA\;X\;Y$ \\
\> \> $(\acomp):: \tA\;X\;Y \to \tA\;Y\;Z \to \tA\;X\;Z$ \\
\> \> $\first:: \tA\;X\;Y \to \tA\;(X,Z)\;(Y,Z)$ \\
\end{tabbing}
satisfying the following equations: 
\begin{center}
\begin{tabular}{|crcl|}
\hline
 (1) & $\arr\; \id \acomp f $ &=& $ f$ \\ 
 (2) & $f \acomp \arr\; \id $ &=& $ f$ \\ 
 (3) & $(f \acomp g) \acomp h $ &=& $ f \acomp (g \acomp h)$ \\ 
 (4) & $\arr\;(w.v) $ &=& $ \arr\; v \acomp \arr\; w$ \\ 
 (5) & $\first\; (\arr\; v) $ &=& $ \arr\; (v\times\id)$ \\
 (6) & $\first\;(f \acomp g) $ &=& $ \first\; f \acomp \first\; g$ \\ 
 (7) & $\first\; f \acomp \arr\; (\id\times v) $ &=& 
    $ \arr\; (\id\times v) \acomp \first\; f$ \\
 (8) & $\first\; f \acomp \arr\; \fst $ &=& $ \arr\; \fst \acomp f$ \\
 (9) & $\;\;\first\; (\first\; f) \acomp \arr\; \assoc $ &=& 
    $ \arr\; \assoc \acomp \first\; f$ \\
\hline
\end{tabular}
\end{center}
where the functions $(\times)$, $\fst$ and $\assoc$ are defined as: 
$$\begin{array}{lll}
(\times) :: & (X\to X')\to(Y\to Y')\to (X,Y)\to (X',Y') &
(f\times g)(x,y)=(f\;x,g\;y) \\
\fst :: & (X,Y)\to X &
\fst(x,y)=x \\
(\assoc) :: & ((X,Y),Z)\to (X,(Y,Z)) &
\assoc((x,y),z) = (x,(y,z)) \\
\end{array}$$
\end{defi}


\subsection{Cartesian effect categories determine Arrows}
\label{subsec:related-cec-arr}

Let $\bV_H$ denote the category of Haskell types and ordinary functions, 
so that the Haskell notation $\mathtt{(X\to Y)}$ represents $\bV_H(X,Y)$, 
made of the Haskell ordinary functions from $X$ to $Y$. 
An arrow $\tA$ contructs a type $\tA\;X\;Y$ for all types $X$ and $Y$.
We slightly modify the definition of Arrows 
by allowing $\mathtt{(X\to Y)}$ to represent $\bV(X,Y)$ for any 
cartesian category $\bV$ 
and by requiring that $\tA\;X\;Y$ is a set rather than a type.
In addition, we use categorical notations instead of Haskell syntax. 

So, from now on, for any cartesian category $\bV$, 
an \emph{Arrow $A$ on $\bV$} associates to each points 
$X$, $Y$ of $\bV$ a set $A(X,Y)$, together with three operations: 
\begin{tabbing}
999 \= 999 \= 999 \kill
\> \> $\arr: \bV(X,Y)\to A(X,Y)$ \\
\> \> $\acomp: A(X,Y) \to A(Y,Z) \to A(X,Z)$ \\
\> \> $\first: A(X,Y) \to A(X\times Z,Y\times Z)$ \\
\end{tabbing}
that satisfy the equations (1)-(9).

Basically, the correspondence between a cartesian effect category $\bC$
extending $\bV$ and an Arrow $A$ on $\bV$ 
identifies $\bC(X,Y)$ with $A(X,Y)$ for all types $X$ and $Y$.
More precisely:

\begin{thm}
\label{thm:related-arr}
Every cartesian effect category $\bC$ extending $\bV$ 
gives rise to an Arrow $A$ on $\bV$, according to the following table:
\textrm{
\begin{center}
\begin{tabular}{|c|c|}
\hline
\boxEff & \boxArr \\
\hline
$\bC(X,Y)$ &
  $A(X,Y)$ \\ 
$\bV(X,Y)\subseteq\bC(X,Y)$ &
  $\arr: \bV(X,Y)\to A(X,Y)$ \\ 
$f\mapsto (g\mapsto g\circ f)$ & 
  $\acomp: A(X,Y) \to A(Y,Z) \to A(X,Z)$ \\
$f\mapsto f\times\id$ & 
  $\first: A(X,Y) \to A(X\times Z,Y\times Z)$  \\ 
\hline
\end{tabular} 
\end{center}
}
\end{thm}

\begin{proof}
The first and second line in the table say that $A(X,Y)$ is 
made of the functions from $X$ to $Y$ in $\bC$
and that $\arr$ is the convertion from pure functions to arbitrary functions. 
The third and fourth lines say that 
$\acomp$ is the (reverse) composition of functions 
and that $\first$ is the semi-product with the identity. 
Let us check that $A$ is an Arrow; 
the following table translates each property (1)-(9) in terms of
cartesian effect categories 
(where $\rho_X:X\times U \to X$ is the projection),
and gives the argument for its proof.
\begin{center}
\begin{tabular}{|crcl|l|}
\hline
 (1) & $f\circ \id $ &=& $f$       & 
  unitarity in $\bC$ \\ 
 (2) & $ \id \circ f $ &=& $f$    & 
  unitarity in $\bC$ \\ 
 (3) & $h\circ (g\circ f) $ &=& $ (h\circ g)\circ f $ & 
  associativity in $\bC$ \\ 
 (4) & $w\circ v$ in $\bV$ &=& $w\circ v$ in $\bC$ &
  $\bV\subseteq\bC$ is a functor \\ 
 (5) & $v\times\id $ in $\bV$ &=& $ v\times\id$ in $\bC$ & 
  non-ambiguity of ``$\times$'' \\ 
 (6) & $(g\circ f)\times\id $ &=& $ (g\times\id)\circ (f\times\id)$ &
  proposition~\ref{prop:effect-comp} \\
 (7) & $(\id\times v) \circ (f\times\id) $ &=& $ (f\times\id) \circ (\id\times v)$ &
  proposition~\ref{prop:effect-comp} \\
 (8) & $\rho\circ (f\times\id_U) $ &=& $ f\circ \rho$ &
  definition~\ref{defi:effect-arr-prod} \\
 (9) & $\alpha^{-1}\circ ((f\times\id)\times\id) $ &=& 
    $ (f\times\id)\circ\alpha^{-1}$ & 
  proposition~\ref{prop:effect-assoc} \\
\hline
\end{tabular}
\end{center}
\end{proof}
\noindent The translation of the Arrow combinators follows easily,
using $\tuple{f,g}_l=(f\ltimes g) \circ \tuple{\id,\id}$ as in section~\ref{subsec:effect-seq}:
\begin{center}
\begin{tabular}{|c|c|}
\hline
\boxEff & \boxArr \\
\hline
$ (id \times f) = \gamma \circ (f \times id)  \circ \gamma  $ &
  $\scond\;f  = \arr\;\swap \acomp \first\;f \acomp \arr\;\swap $ \\
$f \ltimes g  =  (\id\times g) \circ (f\times\id)$ &
  $f \altimes g = \first\;f \acomp \scond\;g $ \\ 
$\tuple{f,g}_l = (f\ltimes g) \circ \tuple{\id,\id}$ &
  $f \atuple g  = \arr(\lambda b \rightarrow (b,b)) \acomp (f \altimes g) $ \\
\hline
\end{tabular} 
\end{center}
For instance, in \cite{Hughes00}, the author states that $\atuple$ 
is not a categorical product since in general 
$(f \atuple g) \acomp \arr\;\fst$ is different from $f$.
We can state this more precisely in the effect category,
where $(f \atuple g) \acomp \arr\;\fst$ corresponds to 
$q_1 \circ \tuple{f,g}_l $.
Indeed, according to corollary~\ref{cor:seq-pair}: 
  $$ q_1 \circ \tuple{f,g}_l  \eqwe  f \;. $$


\section{Examples}
\label{sec:exam}

Here are presented some examples of strict cartesian effect categories.  
Several versions are given, some of them rely on monads.


\subsection{Partiality}
\label{subsec:exam-partial}

Let $\bV=\Set$ be the category of sets and maps,
and $\bC=\Part$ the category of sets and partial maps,
so that $\bV$ is a wide subcategory of $\bC$.
Let $\eqwe$ denote the usual ordering on partial maps:
$f\eqwe g$ if and only if 
$\cD(f)\subseteq\cD(g)$ (where $\cD$ denotes 
the domain of definition) and $f(x)=g(x)$ for all $x\in\cD(f)$. 
The restriction of $\eqwe$ to $\bV$ is the equality of total maps.
Clearly $\eqwe$ is not symmetric, but it satisfies all the other properties 
of a congruence, in particular the 
replacement property with respect to all maps.
So, $\eqwe$ is a semi-congruence (which satisfies replacement), 
that makes $\bC$ a strict effect category extending $\bV$.
\emph{Warning:}  
usually the notations are $v:X\to Y$ for a total map and 
$f:X\parto Y$ for a partial map,
but here we use respectively 
$v:X\sto Y$ (total) and $f:X\to Y$ (partial).

Let us define the pair $\tuple{f,v}$ of a partial map 
$f:X\to Y_1$ and a total map $v:X\sto Y_2$ 
as the partial map $\tuple{f,v}:X\to Y_1\times Y_2$
with the same domain of definition as $f$
and such that $\tuple{f,v}(x)=\tuple{f(x),v(x)}$ for all $x\in\cD(f)$.
It is easy to check that we get a cartesian effect category.
For illustrating the semi-product $f\times v$, there are two cases:
either $f(x_1)$ is defined, or not, in which case we note $f(x_1)=\bot$.
We use the traditional notation 
$\xymatrix@=1.5pc{ x \ar@{|->}[r]^{f} & y \\ }$ when $y=f(x)$
and its analog 
$\xymatrix@=1.5pc{ x \ar@{|~>}[r]^{v} & y \\ }$ when $y=v(x)$
and $v$ is pure.
$$ 
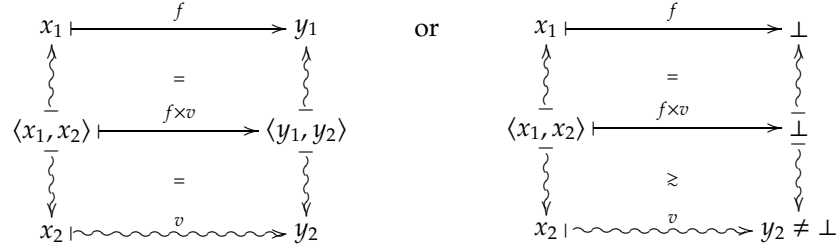
\xymatrix@C=5pc{ 
  x_1 \ar@{|->}[r]^{f} & y_1  \\
  \tuple{x_1,x_2} \ar@{|~>}[u] \ar@{|~>}[d] \ar@{|->}[r]^{f\times v} 
    \ar@{}[rd]|{=} \ar@{}[ru]|{=} & 
    \tuple{y_1,y_2} \ar@{|~>}[u] \ar@{|~>}[d]  \\
  x_2 \ar@{|~>}[r]^{v} & y_2 \\
  }  \qquad \mbox{ or }\qquad
   \xymatrix@C=5pc{ 
  x_1 \ar@{|->}[r]^{f} & \bot  \\
  \tuple{x_1,x_2} \ar@{|~>}[u] \ar@{|~>}[d] \ar@{|->}[r]^{f\times v} 
    \ar@{}[rd]|{\weeq} \ar@{}[ru]|{=} & 
    \bot \ar@{|~>}[u] \ar@{|~>}[d]  \\
  x_2 \ar@{|~>}[r]^{v} & y_2\ne \bot \\
  } 
$$

It can be noted that, in the previous example, $\bC$ is a 2-category, 
with a 2-cell from $f$ to $g$ if and only if $f\eqwe g$.
More generally, let $\bC$ be a 2-category and $\bV$ a sub-2-category 
where the unique 2-cells are the identities.
Then by defining $f\eqwe g$ whenever there is a 2-cell from $f$ to $g$,
we get a strict effect category.
In such effect categories, 
the replacement property holds with respect to all functions in $\bC$,
but the semi-congruence is usually not symmetric.

Let us come back to the partiality example, from the slightly different point of view
of the \emph{Maybe monad}.
First, let us present this point of view in a naive way,
without monads.
Let $U=\{\bot\}$ be a singleton, 
let  ``$+$'' denote the disjoint union of sets,
and for each set $X$ let $GX=X+U$ and let $\eta_X:X\to GX$ be the inclusion.
Each partial map $f$ from $X$ to $Y$ 
can be extended as a total map $Gf$ from $X$ to $GY$,
such that $Gf(x)=f(x)$ for $x\in\cD(f)$ and $Gf(x)=\bot$ otherwise.
This defines a bijection between the partial maps from $X$ to $Y$ 
and the total maps from $X$ to $GY$.
Let $\bC$ be the category such that 
its points are the sets, and a function $X\to Y$ in $\bC$
is a function $X\to GY$ in $\Set$;
we say that $X\to Y$ in $\bC$ \emph{stands for} $X\to GY$ in $\Set$.
Let $J:\Set\to\bC$ be the functor that is the identity on points and 
associates to each map $v_0:X\to Y$ the map $\eta_Y\circ v_0$.
Let $\bV=J(\Set)$.
Then $\bV$ is a wide subcategory of $\bC$.
For all $f,g:X\to Y$ in $\bC$, that stand for $f,g:X\to GY$ in $\Set$,
let:
  $$ f\eqwe g \iff \forall x\in X\; (f(x)\ne\bot\Rightarrow (g(x)\ne\bot \wedge g(x)=f(x))\;.$$
This yields a strict effect category $\bC$ extending $\bV$, 
with the semi-congruence $\eqwe$,
and as above 
the replacement property holds with respect to all functions in $\bC$ 
but $\eqwe$ is not symmetric.
Let $f:X\to Y_1$ in $\bC$ and $v:X\to Y_2$ in $\bV$,
they stand respectively for $f:X\to GY_1$ 
and $v=\eta_{Y_2}\circ v_0$ with $v_0:X\to Y_2$. 
Then, in $\Set$, the pair $\tuple{f,v_0}:X\to GY_1 \times Y_2$ 
can be composed with:
  $$ t:GY_1\times Y_2=(Y_1+U)\times Y_2\to (Y_1\times Y_2)+U =G(Y_1\times Y_2)\;,$$
that maps $\tuple{y_1,y_2}$ to itself and $\tuple{\bot,y_2}$ to $\bot$.
Now, let $\tuple{f,v}:X\to Y_1\times Y_2$ in $\bC$
stand for $\tuple{f,v}=t\circ \tuple{f,v_0}:X\to G(Y_1\times Y_2)$ in $\Set$.
Then $\tuple{f,v}$ is a semi-product, 
so that $\bC$ is a cartesian effect category.
The diagrams for illustrating the semi-product $f\times v$
are the same as above.

This point of view can also be presented using the the \emph{Maybe monad} 
for managing failures, as follows.
We have defined a functor $G:\Part\to\Set$, 
that is a right adjoint to the inclusion functor $I:\Set\subseteq\Part$.
The corresponding monad has endofunctor $M=GI$ on $\Set$,
the category $\bC$ is the Kleisli category of $M$, 
and $J:\Set\to\bC$ is the canonical functor associated to the monad.
In addition, this monad $M$ is \emph{strong}, 
and $t$ is the $(Y_1,Y_2)$ component of the \emph{strength} of $M$. 
But the definition of the semi-congruence $\eqwe$,
as above, is not part of the usual framework of monads.


\subsection{State}
\label{subsec:exam-state}

Let $\bV_0$ be a cartesian category, with a distinguished 
point $S$ for ``the type of states''; for all $X$, 
let $\pi_X:S\times X\to X$ denotes the projection.
Let $\bC$ be the category with the same points as $\bV_0$
and with a function $f:X\to Y$ for each function 
$f:S\times X\to S\times Y$ in $\bV_0$;
we say that $f:X\to Y$ in $\bC$ \emph{stands for} 
$f:S\times X\to S\times Y$ in $\bV_0$.
Let  $J:\bV_0\to\bC$ be the identity-on-points functor
which maps each $v_0:X\to Y$ in $\bV_0$ to
the function $J(v_0):X\to Y$ in $\bC$ that stands for 
$\id_S\times v_0:S\times X\to S\times Y$ in $\bV_0$.
Let $\bV=J(\bV_0)$, it is a wide subcategory of $\bC$.
For all $f,g:X\to Y$ in $\bC$, let:
  $$ f\eqwe g \iff \pi_Y\circ g = \pi_Y\circ f \;. $$
We get a strict effect category, where the semi-congruence $\eqwe$ 
is symmetric, but does not satisfy the replacement property
with respect to all functions in $\bC$.
The semi-product of $f:X\to Y_1$ and $v:X\sto Y_2$ is defined as follows.
Since $f:S\times X\to S\times Y_1$ in $\bV_0$
and $v=\id_S\times v_0$ for some $v_0:X\to Y$ in $\bV_0$,
the pair $\tuple{f,v_0\circ\pi_X}:S\times X\to (S\times Y_1)\times Y_2$ 
exists in $\bV_0$. By composing it with the isomorphism 
$(S\times Y_1)\times Y_2 \to S\times (Y_1\times Y_2)$
we get $\tuple{f,v}:S\times X\to S\times (Y_1\times Y_2)$ in $\bV_0$,
i.e., $\tuple{f,v}:X\to Y_1\times Y_2$ in $\bC$.
It is easy to check that this defines a semi-product, 
so that $\bC$ is a cartesian effect category,
where the characteristic property of the semi-product $f\times v$
can be illustrated as follows:
$$ \xymatrix@C=5pc{ 
  (s,x_1) \ar@{|->}[r]^{f} & (s',y_1)  \\
  (s,x_1,x_2) \ar@{|~>}[u] \ar@{|~>}[d] \ar@{|->}[r]^{f\times v} 
    \ar@{}[rd]|{\weeq} \ar@{}[ru]|{=} & 
    (s',y_1,y_2) \ar@{|~>}[u] \ar@{|~>}[d]  \\
  (s,x_2) \ar@{|~>}[r]^{v} & (s,y_2)\ne(s',y_2) 
    \ar@{|~>}[r]^{\quad\pi_Y} & y_2 \\
  }  $$

The example above can be curried, thus recovering the \emph{State monad}.
A motivation for the introduction of Freyd categories 
in \cite{PowerRobinson97} is the possibility 
of dealing with state in a linear way, as above,
rather than in the exponential way provided by the State monad.
Now $\bV_0$ is still a cartesian category with a distinguished 
point $S$, the ``type of states'',
and in addition $\bV_0$ has exponentials $(S\times X)^S$ for each $X$. 
Then the endofunctor $M(X)=(S\times X)^S$ defines the State monad on $\bV_0$, 
with composition defined as usual. 
It is well-known that $M$ is a strong monad, with strength 
$t_{Y_1,Y_2}=(S\times Y_1)^S\times Y_2 \to (S\times Y_1\times Y_2)^S$
obtained from 
$\app_{S\times Y_1}\times\id_{Y_2}:
S\times (S\times Y_1)^S\times Y_2 \to S\times Y_1\times Y_2$,
where ``$\app$'' denotes the application function. 
Hence, from $f:X\to M(Y_1)$ and $v_0:X\to Y_2$ in $\bV_0$, 
we can build $\tuple{f,v}=t_{Y_1,Y_2}\circ\tuple{f,v_0}:X\to M(Y_1\times Y_2)$.
Let $\bC$ be the Kleisli category of the monad $M$,
let $J:\bV_0\to\bC$ be the canonical functor associated to the monad,
and let $\bV=J(\bV_0)$,
then $\bV$ is a wide subcategory of $\bC$. 
A function $f:X\to Y$ in $\bC$ stands for a function $f:X\to (S\times Y)^S$ in $\bV_0$.
Now, in addition to the usual framework of monads, 
for all $f,g:X\to Y$ in $\bC$, i.e., $f,g:X\to (S\times Y)^S$ in $\bV_0$,
let:
  $$ f\eqwe g \iff  {\pi_Y}^S \circ g = {\pi_Y}^S \circ f \;, $$
where ${\pi_Y}^S:(S\times Y)^S\to Y^S$ associates to each 
map $m:S\to S\times Y$ the map $\pi_Y\times m:S\to Y$.
The relation $\eqwe$ defines a semi-conguence on $\bC$,
and $\tuple{f,v}$ is a semi-product, 
so that $\bC$ is a cartesian effect category.
The characteristic property of the semi-product $f\times v$
can be illustrated as follows:
$$ \xymatrix@C=5pc{ 
  x_1 \ar@{|->}[r]^{f} & (s\mapsto(s',y_1))  \\
  (x_1,x_2) \ar@{|~>}[u] \ar@{|~>}[d] \ar@{|->}[r]^{f\times v} 
    \ar@{}[rd]|{\weeq} \ar@{}[ru]|{=} & 
    (s\mapsto(s',y_1,y_2)) \ar@{|~>}[u] \ar@{|~>}[d]  \\
  x_2 \ar@{|~>}[r]^{v} & (s\mapsto(s,y_2))\ne(s\mapsto(s',y_2)) 
    \ar@{|~>}[r]^{\qquad\pi_Y^{S}} & (s\mapsto y_2) \\
  }  $$


\section{Conclusion}

We have presented a new categorical framework,
called a \emph{cartesian effect category}, 
for dealing with the issue of multiple arguments 
in programming languages with computational effects. 
The major new feature in cartesian effect categories is the 
introduction of a \emph{semi-congruence},
which allows to define \emph{semi-products} 
and to prove their properties 
by decorating the usual definitions, properties 
and proofs about products in a category.
Forthcoming work should study the nesting of several effects. 

In order to deal with other issues related to effects, 
we believe that the idea of \emph{decorations} in logic
can be more widely used. This is the case for 
dealing with exceptions \cite{B-DuvRey05}
(note that a previous attempt to define decorated products
can be found in \cite{F-DuvRey04}).
The framework of decorations might be used 
for generalizing this work in the direction of 
closed Freyd categories \cite{PowerThielecke99}.
or traced premonoidal categories \cite{BentonHyland03}.
Moreover, with one additional level of abstraction,
decorations can be obtained from \emph{morphisms between logics}, 
in the context of \emph{diagrammatic logics} \cite{F-DuvLai02,A-Duv03}.


\bibliographystyle{plain}
\bibliography{prod}


\appendix


\section{Proofs in cartesian effect categories}
\label{app:proof}

Here are proofs for some results in section~\ref{subsec:weak-prod},
called \emph{basic} proofs, 
followed by their \emph{decorated} versions for the corresponding 
results in section~\ref{subsec:effect-semi}.
All basic proofs are straightforward.
All proofs are presented in a formalized way: 
each property is preceded by its label and followed by its proof.
For the basic proofs, the properties of the congruence are denoted 
$\tr$, $\sym$, $\subst$, $\repl$, for respectively 
transitivity, symmetry, substitution, replacement. 
For the decorated proofs, 
the properties of the congruence and the semi-congruence are still denoted 
$\tr$, $\sym$, $\subst$, $\repl$, with subscript either $\eqst$ or $\eqwe$.
It should be reminded that $\sym_{\eqwe}$ does \emph{not} hold,
and that $\repl_{\eqwe}$ is allowed \emph{only} with respect to a pure function:
if $g_1\eqwe g_2:Y\to Z$ and $v:Z\sto W$ 
then $v\circ g_1 \eqwe v\circ g_2:Y\to W$.
In addition, $\comp$ means compatibiblity of $\eqwe$ with $\eqst$,
which means that 
if either $f_1\eqst f_2\eqwe f_3$ or $f_1\eqwe f_2\eqst f_3$ 
then $f_1\eqwe f_3$.
In decorated proofs, 
``like \emph{basic}'' means that this part of the proof 
is exactly the same as in the basic proof.
Proofs of propositions~\ref{prop:weak-assoc}, \ref{prop:effect-assoc}(associativity) 
and~\ref{prop:weak-seq}, \ref{prop:effect-seq} (parallelism) are left to the reader.

\begin{proof}[Proof of proposition~\ref{prop:weak-equiv} (congruence)]\
\\ \begin{tabular}{rll}
1.$\quad$  & When $X_1=X_2$ & \\
$(a_1)$ & $q_1\circ \tuple{f_1,f_2} \equiv f_1$  & \\ 
$(b_1)$ & $f_1 \equiv f'_1$ & \\ 
$(c_1)$ & $q_1\circ \tuple{f_1,f_2} \equiv f'_1$ & $(a_1)$, $(b_1)$, $\tr$ \\
$(c_2)$ & $q_2\circ \tuple{f_1,f_2} \equiv f'_2$ & like $(c_1)$ \\
$(d)$ & $\tuple{f_1,f_2} \equiv \tuple{f'_1,f'_2}$ & $(c_1)$, $(c_2)$ \\
2.$\quad$ & In all cases & \\ 
$(e_1)$ & $f_1\equiv f'_1$ & \\ 
$(f_1)$ & $f_1\circ p_1\equiv f'_1\circ p_1$ & $(e_1)$, $\subst$ \\
$(f_2)$ & $f_2\circ p_2\equiv f'_2\circ p_2$ & like $(f_1)$ \\
$(g)$ & $\tuple{f_1\circ p_1,f_2\circ p_2} \equiv 
\tuple{f'_1\circ p_1,f'_2\circ p_2}$ & $(f_1)$, $(f_2)$, $(1)$ \\
\end{tabular}\\ 
\end{proof}

\begin{proof}[Proof of proposition~\ref{prop:effect-equiv} (congruence)]\
\\ \begin{tabular}{rll}
1.$\quad$  & When $X_1=X_2$ & \\
$(c_1)$ & $q_1\circ \tuple{f_1,f_2} \eqst f'_1$ & like \emph{basic} \\
$(a_2)$ & $q_2\circ \tuple{f_1,f_2} \eqwe f_2$  & \\ 
$(b_2)$ & $f_2 \eqst f'_2$ & \\ 
$(c_2)$ & $q_2\circ \tuple{f_1,f_2} \eqwe f'_2$ & $(a_2)$, $(b_2)$, $\comp$ \\
$(d)$ & $\tuple{f_1,f_2} \eqst \tuple{f'_1,f'_2}$ & $(c_1)$, $(c_2)$ \\
2.$\quad$ & In all cases & \\ 
$(g)$ & $\tuple{f_1\circ p_1,f_2\circ p_2} \eqst 
\tuple{f'_1\circ p_1,f'_2\circ p_2}$ & like \emph{basic} \\
\end{tabular}\\ 
\end{proof}

\begin{proof}[Proof of proposition~\ref{prop:weak-comp} (composition)]\hfill 
\\ The three left handsides can be illustrated as follows:
  $$ 
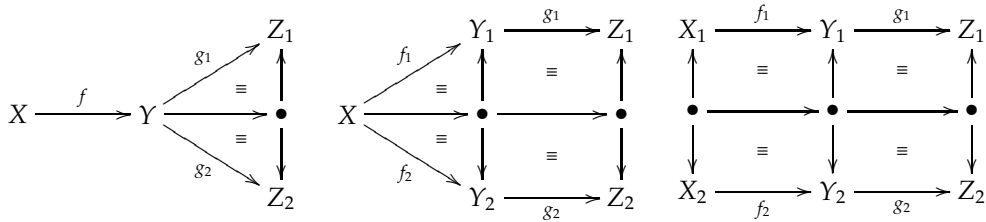
\xymatrix@C=3pc@R=1.5pc{ 
  & & Z_1 \\
  X \ar[r]^{f} & Y \ar[ru]^{g_1} \ar[rd]_{g_2} \ar[r] &
    \bullet \ar[u] \ar[d] 
    \ar@{}[ld]|(.3){\equiv} \ar@{}[lu]|(.3){\equiv} \\
  & & Z_2 \\
  } \quad
  \xymatrix@C=3pc@R=1.5pc{ 
  & Y_1 \ar[r]^{g_1} & Z_1 \\
  X \ar[ru]^{f_1} \ar[rd]_{f_2} \ar[r] &
    \bullet \ar[u] \ar[d] \ar[r] 
    \ar@{}[ld]|(.3){\equiv} \ar@{}[lu]|(.3){\equiv} &
    \bullet \ar[u] \ar[d] 
    \ar@{}[ld]|{\equiv} \ar@{}[lu]|{\equiv} \\
  & Y_2 \ar[r]_{g_2} & Z_2 \\
  } \quad
  \xymatrix@C=3pc@R=1.5pc{ 
  X_1\ar[r]^{f_1} & Y_1 \ar[r]^{g_1} & Z_1 \\
  \bullet \ar[u] \ar[d] \ar[r] &
    \bullet \ar[u] \ar[d] \ar[r] 
    \ar@{}[ld]|{\equiv} \ar@{}[lu]|{\equiv} &
    \bullet \ar[u] \ar[d] 
    \ar@{}[ld]|{\equiv} \ar@{}[lu]|{\equiv} \\
  X_2 \ar[r]_{f_2} & Y_2  \ar[r]_{g_2} & Z_2 \\
  } $$
\begin{tabular}{rll}
1.$\quad$ & When $f_1=f_2(=f)$ & \\
$(a_1)$ & $r_1\circ \tuple{g_1,g_2} \equiv g_1$ & \\ 
$(b_1)$ & $r_1\circ\tuple{g_1,g_2}\circ f\equiv g_1\circ f $ & $(a_1)$, $\subst$ \\
$(b_2)$ & $r_2\circ \tuple{g_1,g_2} \circ f \equiv g_2 \circ f $ & like $(b_1)$ \\
$(c)$ & $\tuple{g_1,g_2} \circ f \equiv \tuple{g_1\circ f,g_2\circ f}$ & $(b_1)$, $(b_2)$ \\
2.$\quad$ & When $X_1=X_2$ & \\ 
$(d)$ & $(g_1\times g_2) \circ \tuple{f_1,f_2} 
\equiv \tuple{g_1\circ q_1 \circ \tuple{f_1,f_2},
g_2\circ q_2 \circ \tuple{f_1,f_2}}$ & $(1)$ \\
$(e_1)$ & $q_1 \circ \tuple{f_1,f_2} \equiv f_1$ & \\
$(f_1)$ & $g_1\circ q_1\circ\tuple{f_1,f_2} \equiv g_1\circ f_1$ & $\repl$ \\
$(f_2)$ & $g_2\circ q_2 \circ \tuple{f_1,f_2}\equiv g_2\circ f_2$ & like $(f_1)$ \\
$(g)$ & $\tuple{g_1\circ q_1 \circ \tuple{f_1,f_2},
g_2\circ q_2 \circ \tuple{f_1,f_2}} \equiv 
\tuple{g_1\circ f_1,g_2\circ f_2}$ & $(f_1)$, $(f_2)$, prop.~\ref{prop:weak-equiv} \\
$(h)$ & $(g_1\times g_2) \circ \tuple{f_1,f_2} \equiv 
\tuple{g_1\circ f_1,g_2\circ f_2}$ & $(d)$, $(g)$, $\tr$ \\
3.$\quad$ & In all cases & \\
$(k)$ & $ (g_1 \times g_2) \circ \tuple{f_1\circ p_1,f_2\circ p_2} \equiv 
      \tuple{g_1\circ f_1\circ p_1,g_2\circ f_2\circ p_2}$ & $(2)$ \\
\end{tabular}\\
\end{proof}

\begin{proof}[Proof of proposition~\ref{prop:effect-comp} (composition)]\hfill 
\\ The three left handsides can be illustrated as follows:
  $$ 
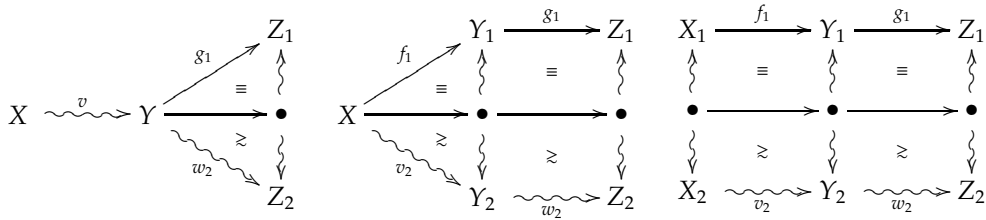
\xymatrix@C=3pc@R=1.5pc{ 
  & & Z_1 \\
  X \ar@{~>}[r]^{v} & Y \ar[ru]^{g_1} \ar@{~>}[rd]_{w_2} \ar[r] &
    \bullet \ar@{~>}[u] \ar@{~>}[d] 
    \ar@{}[ld]|(.3){\weeq} \ar@{}[lu]|(.3){\eqst} \\
  & & Z_2 \\
  } \quad
  \xymatrix@C=3pc@R=1.5pc{ 
  & Y_1 \ar[r]^{g_1} & Z_1 \\
  X \ar[ru]^{f_1} \ar@{~>}[rd]_{v_2} \ar[r] &
    \bullet \ar@{~>}[u] \ar@{~>}[d] \ar[r] 
    \ar@{}[ld]|(.3){\weeq} \ar@{}[lu]|(.3){\eqst} &
    \bullet \ar@{~>}[u] \ar@{~>}[d] 
    \ar@{}[ld]|{\weeq} \ar@{}[lu]|{\eqst} \\
  & Y_2 \ar@{~>}[r]_{w_2} & Z_2 \\
  } \quad
  \xymatrix@C=3pc@R=1.5pc{ 
  X_1\ar[r]^{f_1} & Y_1 \ar[r]^{g_1} & Z_1 \\
  \bullet \ar@{~>}[u] \ar@{~>}[d] \ar[r] &
    \bullet \ar@{~>}[u] \ar@{~>}[d] \ar[r] 
    \ar@{}[ld]|{\weeq} \ar@{}[lu]|{\eqst} &
    \bullet \ar@{~>}[u] \ar@{~>}[d] 
    \ar@{}[ld]|{\weeq} \ar@{}[lu]|{\eqst} \\
  X_2 \ar@{~>}[r]_{v_2} & Y_2  \ar@{~>}[r]_{w_2} & Z_2 \\
  } $$
\begin{tabular}{rll}
1.$\quad$ & When $f_1=v_2(=v)$ & \\
$(b_1)$ & $r_1\circ\tuple{g_1,w_2}\circ v\eqst g_1\circ v $ & like \emph{basic} \\
$(a_2)$ & $r_2\circ \tuple{g_1,w_2} \eqwe w_2 $ & \\ 
$(b_2)$ & $r_2\circ\tuple{g_1,w_2}\circ v\eqwe w_2\circ v $ & $(a_1)$, $\subst_{\eqwe}$ \\
$(c)$ & $\tuple{g_1,w_2} \circ v \eqst \tuple{g_1\circ v,w_2\circ v}$ & $(b_1)$, $(b_2)$ \\
2.$\quad$ & When $X_1=X_2$ & \\ 
$(d)$ & $(g_1\times w_2) \circ \tuple{f_1,v_2} 
\eqst \tuple{g_1\circ q_1 \circ \tuple{f_1,v_2},
w_2\circ q_2 \circ \tuple{f_1,v_2}}$ & $(1)$ \\
$(f_1)$ & $g_1\circ q_1\circ\tuple{f_1,v_2} \eqst g_1\circ f_1$ & like \emph{basic} \\
$(e_2)$ & $q_2 \circ \tuple{f_1,v_2} \eqwe v_2$ & \\
$(f_2)$ & $w_2\circ q_2\circ\tuple{f_1,v_2} \eqwe w_2\circ v_2$ & $\repl_{\eqwe}$ ($w_2$ is pure)\\
$(g)$ & $\tuple{g_1\circ q_1 \circ \tuple{f_1,v_2},
w_2\circ q_2 \circ \tuple{f_1,v_2}} \eqst 
\tuple{g_1\circ f_1,w_2\circ v_2}$ & $(f_1)$, $(f_2)$, prop.~\ref{prop:weak-equiv} \\
$(h)$ & $(g_1\times w_2) \circ \tuple{f_1,v_2} \eqst
\tuple{g_1\circ f_1,w_2\circ v_2}$ & $(d)$, $(g)$, $\tr_{\eqst}$ \\
3.$\quad$ & In all cases & \\
$(k)$ & $ (g_1 \times w_2) \circ \tuple{f_1\circ p_1,v_2\circ p_2} \eqst
      \tuple{g_1\circ f_1\circ p_1,w_2\circ v_2\circ p_2}$ & $(2)$ \\
\end{tabular}\\
\end{proof}

\begin{proof}[Proof of proposition~\ref{prop:weak-swap} (swap)]\hfill 
\\ The two left handsides can be illustrated as follows:
  $$ 
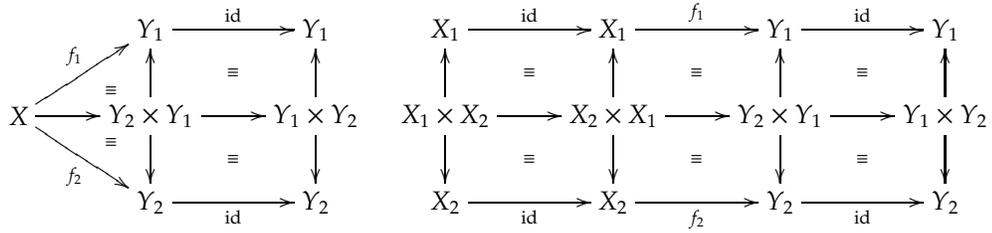
\xymatrix@C=2pc@R=1.5pc{ 
  & Y_1 \ar[r]^{\id} & Y_1 \\
  X \ar[ur]^{f_1} \ar[dr]_{f_2} \ar[r] & 
  Y_2\times Y_1 \ar[u]\ar[d] \ar[r] \ar@{}[ld]|(.3){\equiv} \ar@{}[lu]|(.3){\equiv} & 
  Y_1\times Y_2 \ar[u]\ar[d] \ar@{}[ld]|{\equiv} \ar@{}[lu]|{\equiv} \\
   & Y_2 \ar[r]_{\id} & Y_2 \\
  }   \quad \xymatrix@C=2pc@R=1.5pc{ 
  X_1 \ar[r]^{\id} & X_1 \ar[r]^{f_1} & Y_1 \ar[r]^{\id} & Y_1 \\
  X_1\times X_2 \ar[u]\ar[d] \ar[r] & 
  X_2\times X_1 \ar[u]\ar[d] \ar[r] \ar@{}[ld]|{\equiv} \ar@{}[lu]|{\equiv} & 
  Y_2\times Y_1 \ar[u]\ar[d] \ar[r] \ar@{}[ld]|{\equiv} \ar@{}[lu]|{\equiv} & 
  Y_1\times Y_2 \ar[u]\ar[d] \ar@{}[ld]|{\equiv} \ar@{}[lu]|{\equiv} \\
  X_2 \ar[r]_{\id} & X_2 \ar[r]_{f_2} & Y_2 \ar[r]_{\id} & Y_2 \\
  } 
$$
\begin{tabular}{rll}
1.$\quad$ & When $X_1=X_2$ & \\
$(a_1)$ & $q_1\circ\gamma_Y \equiv q'_1$ & \\
$(b_1)$ & $q_1\circ\gamma_Y \circ \tuple{f_2,f_1} 
\equiv q'_1 \circ \tuple{f_2,f_1} $ & $(a_1)$, $\subst$ \\
$(c_1)$ & $q'_1 \circ \tuple{f_2,f_1} \equiv f_1$ & \\
$(d_1)$ & $q_1\circ\gamma_Y \circ \tuple{f_2,f_1} \equiv f_1$ &
    $(b_1)$, $(c_1)$, $\tr$ \\
$(d_2)$ & $q_2\circ\gamma_Y \circ \tuple{f_2,f_1} \equiv f_2 $ & like $(d_1)$ \\
$(e)$ & $\gamma_Y \circ \tuple{f_2,f_1} \equiv \tuple{f_1,f_2}$ & $(d_1)$, $(d_2)$ \\
2.$\quad$ & In all cases & \\ 
$(f)$ & $\tuple{f_2\circ p'_2, f_1\circ p'_1} \circ \gamma_X^{-1}   
  \equiv \tuple{f_2\circ p'_2\circ \gamma_X^{-1}, f_1\circ p'_1\circ \gamma_X^{-1}} $
 & prop.~\ref{prop:weak-comp}, $\sym$ \\
$(g)$ & $\gamma_Y \circ \tuple{f_2\circ p'_2, f_1\circ p'_1} \circ \gamma_X^{-1}   
  \equiv \gamma_Y \circ \tuple{f_2\circ p'_2\circ \gamma_X^{-1}, f_1\circ p'_1\circ \gamma_X^{-1}} $
 &$\repl$ \\
$(h)$  & $\gamma_Y \circ \tuple{f_2\circ p'_2\circ \gamma_X^{-1}, f_1\circ p'_1\circ \gamma_X^{-1}}
  \equiv \tuple{f_1\circ p'_1\circ \gamma_X^{-1},f_2\circ p'_2\circ \gamma_X^{-1}}$
 & $(1)$ \\
$(i_1)$ & $p'_1\circ \gamma_X^{-1} \equiv p_1$ & \\
$(j_1)$ & $f_1\circ p'_1\circ \gamma_X^{-1} \equiv f_1\circ p_1$ & $(i_1)$, $\repl$ \\
$(j_2)$ & $f_2\circ p'_2\circ \gamma_X^{-1} \equiv f_2\circ p_2$ & like $(j_1)$ \\
$(k)$ &$ \tuple{f_1\circ p'_1\circ \gamma_X^{-1},f_2\circ p'_2\circ \gamma_X^{-1}}
\equiv \tuple{f_1\circ p_1,f_2\circ p_2}$ & $(j_1)$, $(j_2)$, prop.~\ref{prop:weak-equiv} \\
$(l)$ & $\gamma_Y \circ \tuple{f_2\circ p'_2, f_1\circ p'_1} \circ \gamma_X^{-1}
\equiv \tuple{f_1\circ p_1,f_2\circ p_2}$ & $(g)$, $(h)$, $(k)$, $\tr$ \\
\end{tabular}\\
\end{proof}

\begin{proof}[Proof of proposition~\ref{prop:effect-swap} (swap)]\hfill 
\\ The two left handsides can be illustrated as follows:
  $$ \xymatrix@C=2pc@R=1.5pc{ 
  & Y_1 \ar@{~>}[r]^{\id} & Y_1 \\
  X \ar[ur]^{f_1} \ar@{~>}[dr]_{v_2} \ar[r] & 
  Y_2\times Y_1 \ar@{~>}[u]\ar@{~>}[d] \ar@{~>}[r] 
     \ar@{}[ld]|(.3){\weeq} \ar@{}[lu]|(.3){\eqst} & 
  Y_1\times Y_2 \ar@{~>}[u]\ar@{~>}[d] \ar@{}[ld]|{\eqst} \ar@{}[lu]|{\eqst} \\
   & Y_2 \ar@{~>}[r]_{\id} & Y_2 \\
  }   \quad \xymatrix@C=2pc@R=1.5pc{ 
  X_1 \ar@{~>}[r]^{\id} & X_1 \ar[r]^{f_1} & Y_1 \ar@{~>}[r]^{\id} & Y_1 \\
  X_1\times X_2 \ar@{~>}[u]\ar@{~>}[d] \ar@{~>}[r] & 
  X_2\times X_1 \ar@{~>}[u]\ar@{~>}[d] \ar[r] \ar@{}[ld]|{\eqst} \ar@{}[lu]|{\eqst} & 
  Y_2\times Y_1 \ar@{~>}[u]\ar@{~>}[d] \ar@{~>}[r] \ar@{}[ld]|{\weeq} \ar@{}[lu]|{\eqst} & 
  Y_1\times Y_2 \ar@{~>}[u]\ar@{~>}[d] \ar@{}[ld]|{\eqst} \ar@{}[lu]|{\eqst} \\
  X_2 \ar@{~>}[r]_{\id} & X_2 \ar@{~>}[r]_{v_2} & Y_2 \ar@{~>}[r]_{\id} & Y_2 \\
  } 
$$
\begin{tabular}{rll}
1.$\quad$ & When $X_1=X_2$ & \\
$(d_1)$ & $q_1\circ\gamma_Y \circ \tuple{v_2,f_1} \eqst f_1$ & like \emph{basic} \\
$(a_2)$ & $q_2\circ\gamma_Y \eqst q'_2$ & \\
$(b_2)$ & $q_2\circ\gamma_Y \circ \tuple{v_2,f_1} 
\eqst q'_2 \circ \tuple{v_2,f_1} $ & $(a_2)$, $\subst_{\eqst}$ \\
$(c_2)$ & $q'_2 \circ \tuple{v_2,f_1} \eqwe v_2$ & \\
$(d_2)$ & $q_2\circ\gamma_Y \circ \tuple{v_2,f_1} \eqwe v_2$ &
    $(b_2)$, $(c_2)$, $\comp$ \\
$(e)$ & $\gamma_Y \circ \tuple{v_2,f_1} \eqst \tuple{f_1,v_2}$ & $(d_1)$, $(d_2)$ \\
2.$\quad$ & In all cases & \\ 
$(l)$ & $\gamma_Y \circ \tuple{f_2\circ p'_2, f_1\circ p'_1} \circ \gamma_X^{-1}
\equiv \tuple{f_1\circ p_1,f_2\circ p_2}$ & like \emph{basic} \\ 
\end{tabular}\\
\end{proof}


\end{document}